\documentclass{paper}

\pdfoutput=1

\usepackage{amsmath}
\usepackage{amsfonts,amsmath}
\usepackage{amssymb}
\usepackage{array}
\usepackage{multirow}
\usepackage{xspace}
\usepackage{float}
\usepackage{subfig}
\usepackage{tikz}
\usetikzlibrary{arrows}
\usepackage{fullpage}

\numberwithin{equation}{section}
\renewcommand\arraystretch{1.5}

\title{Certified Reduced Basis Method for Affinely Parametric Isogeometric Analysis NURBS Approximation}
\author{Denis Devaud \thanks{Universit\'e de Neuch\^atel, Institut de Statistique, Avenue de Bellevaux 51, 2000 Neuch\^atel, Switzerland, \texttt{denis.devaud@unine.ch}} \and Gianluigi Rozza \thanks{SISSA, International School for Advanced Studies, Mathematics Area, mathLab, Via Bonomea 265, 34136 Trieste, Italy, \texttt{grozza@sissa.it}}}

\begin{document}
\maketitle

\section{Introduction and Motivation}
In this work we apply reduced basis methods for parametric PDEs to an isogeometric formulation based on NURBS. The motivation for this work is an integrated and complete work pipeline from CAD to parametrization of domain geometry, then from full order to certified reduced basis solution. IsoGeometric Analysis (IGA) is a growing research theme in scientific computing and computational mechanics, as well as reduced basis methods for parametric PDEs. Their combination enhances the solution of some class of problems, especially the ones characterized by parametrized geometries we introduced in this work.
For a general overview on Reduced Basis (RB) methods we recall \cite{hesthaven2015certified, rozza2011reduced} and on IGA \cite{cottrell2009isogeometric}.
This work wants to demonstrate that it is also possible for some class of problems to deal with affine geometrical parametrization combined with a NURBS IGA formulation. This is what this work brings as original ingredients with respect to other works dealing with reduced order methods and IGA (set in a non-affine formulation, and using a POD \cite{RozzaEncyclopedia} sampling without certification: see for example for potential flows \cite{manzoni2015reduced} and for Stokes flows \cite{Salmoiraghi2016viscous}). In this work we show a certification of accuracy and a complete integration between IGA formulation and parametric certified greedy RB formulation.
Section 2 recalls the abstract setting for parametrized PDEs, Section 3 recalls IGA setting, Section 4 deals with RB formulation, and Section 5 illustrates two numerical examples in heat transfer with different parametrization.

\section{Elliptic Coercive Parametrized Partial Differential Equations}

In what follows, elliptic coercive parametrized partial differential equations are introduced \cite{patera2007reduced,quarteroni2011certified,rozza2008reduced}.
We consider the following problem: given a parameter $\boldsymbol\mu\in\mathcal{D}$, evaluate
\begin{equation}\label{equ:sex}
s(\boldsymbol\mu)=l(u(\boldsymbol\mu)),
\end{equation}
where $u(\boldsymbol\mu)\in X$ is the solution of
\begin{equation}\label{equ:aex}
a(u(\boldsymbol\mu),v;\boldsymbol\mu)=f(v),\qquad\forall v\in X.
\end{equation}
Here $a(\cdot,\cdot;\boldsymbol\mu):X\times X\rightarrow\mathbb{R}$ is a bilinear, continuous and coercive form associated to a parametrized partial differential equation for every $\boldsymbol\mu\in\mathcal{D}$.
The space $X:=X(\Omega)$ is a Hilbert space on the computational domain $\Omega\subset\mathbb{R}^{d}$ endowed with the scalar product $(\cdot,\cdot)_{X}$ for $d=2,3$.
Since second-order partial differential equations for scalar problems are considered, we have $H_{0}^{1}(\Omega)\subset X\subset H^{1}(\Omega)$, where $H^{1}(\Omega):=\left\{v:\Omega\rightarrow X\middle|\ v\in L^{2}(\Omega),\nabla v\in L^{2}(\Omega)^{d}\right\}$ and $H_{0}^{1}(\Omega)$ is the space of functions in $H^{1}(\Omega)$ whose traces vanish on the boundary.
The space $L^{2}(\Omega)$ denotes the set of square integrable functions.
We require moreover that $\Omega$ admits a (multipatches) NURBS representation.
This is explained in more details in the next section.
The functions $f:X\rightarrow\mathbb{R}$ and $l:X\rightarrow\mathbb{R}$ are linear and continuous functionals.
Finally, the set $\mathcal{D}$ denotes the parameter domain and is assumed to be finite-dimensional.
More precisely, we write $\mathcal{D}:=[a_{1},b_{1}]\times\dots\times[a_{P},b_{P}]\subset\mathbb{R}^{P}$ for $a_{i},b_{i}\in\mathbb{R}$, $i=1,\dots,P$.
We consider here both physical and geometrical parameters.
The geometrical case is further investigated in Section \ref{subsec:affine_precond}.
For a sake of simplicity, the so-called compliant case is considered, that is $(i)$ $a$ is symmetric and $(ii)$ $l=f$.

One of the crucial assumptions to apply the reduced basis method is that $a$ admits an affine decomposition with respect to the parameter $\boldsymbol\mu$, that is
\begin{equation}\label{equ:affinea}
a(u,v;\boldsymbol\mu) = \sum_{q=1}^{Q}\Theta^{q}(\boldsymbol\mu)a^{q}(u,v).
\end{equation}
Here $\Theta^{q}:\mathcal{D}\rightarrow\mathbb{R}$ denotes a (smooth) $\boldsymbol\mu$-dependent function and $a^{q}:X\times X\rightarrow\mathbb{R}$  is a $\boldsymbol\mu$-independent bilinear continuous form for $q=1,\dots,Q$.
Since the compliant case is considered, we require moreover that $a^{q}$ is symmetric.
We do not make any further assumption on the coercivity of $a^{q}$.
Note that we have assumed that the right-hand side of equation \eqref{equ:aex} is parameter-independent but in practice $f$ may depend on the parameter $\boldsymbol\mu$.
In that case, we express $f(v;\boldsymbol\mu)$ as a sum of $Q_{f}$ products of $\boldsymbol\mu$-dependent functions and $\boldsymbol\mu$-independent linear continuous forms on $X$.

For the bilinear form $a(\cdot,\cdot;\boldsymbol\mu)$, we define its continuity and coercivity constants for every $\boldsymbol\mu\in\mathcal{D}$ as
\begin{equation}\label{equ:excontconstant}\nonumber
\gamma(\boldsymbol\mu):=\sup_{v\in X}\sup_{w\in X}\frac{a(v,w;\boldsymbol\mu)}{\|v\|_{X}\|w\|_{X}},
\end{equation}
and
\begin{equation}\label{equ:excoerconstant}\nonumber
\alpha(\boldsymbol\mu):=\inf_{v\in X}\frac{a(v,v;\boldsymbol\mu)}{\|v\|_{X}^{2}},
\end{equation}
where $\|\cdot\|_{X}$ is the norm on $X$ induced by the scalar product $(\cdot,\cdot)_{X}$.
Since $a$ is continuous and coercive, there exists $0<\alpha_{0}\leq\gamma_{0}<\infty$ such that $\alpha_{0}\leq\alpha(\boldsymbol\mu)\leq\gamma(\boldsymbol\mu)\leq\gamma_{0}$ for all $\boldsymbol\mu\in\mathcal{D}$.

In the following section, we introduce a NURBS approximation of the problem \eqref{equ:sex}-\eqref{equ:aex}.
Since it is computationally unaffordable to compute such solution for every input parameter, we then consider a RB approximation of it.

\section{Isogeometric Analysis NURBS Approximation}

In this section, we introduce non-uniform rational B-splines (NURBS) approximation for the problem \eqref{equ:sex}-\eqref{equ:aex}.
First, a brief survey of B-splines and NURBS functions is conducted and the proper approximation is introduced \cite{bazilevs2006isogeometric,cottrell2009isogeometric,hughes2005isogeometric}.
We then present in Section \ref{subsec:affine_precond} the affine preconditioning for parameter-dependent domains.
This is a necessary assumption to obtain the affine decomposition \eqref{equ:affinea} which in turn is crucial to perform RB approximation.

\subsection{B-Splines}\label{subsec:bsplines}

The B-splines functions are the basis to define NURBS.
We give a brief introduction to B-splines in what follows.
In the context of isogeometric analysis, the notion of \emph{patches} is very important.
They play the role of subdomains and material properties are assumed to be uniform in each patch.
Unlike standard finite element (FE) analysis, the B-splines and NURBS basis functions are local to patches and not elements.
The FE basis functions map the reference element in the parametric domain to each element in the physical space.
B-splines functions take a patch (a set of elements) in the parameter space and map it to multiple elements in the physical domain.

Let us define a \emph{knot vector} in one dimension as a set of non-decreasing coordinates in the parameter domain denoted $\Xi=\left\{\xi_{1},\dots,\xi_{n+p+1}\right\}$, where $\xi_{i}\in\mathbb{R}$ is called the $i$th \emph{knot}, $i=1,\dots,n+p+1$.
Here, $p$ denotes the polynomial order of the B-splines and $n$ the number of basis functions.
The B-splines are completely defined by the knot vector $\Xi$, the number of basis functions $n$ and their order $p$.
Since this does not affect the construction of B-splines we set by convention $\xi_{1}=0$ and $\xi_{n+p+1}=1$.
Note that repetitions are allowed in the knot vector and are used to control the local regularity across each knot.
A knot vector in which $\xi_{1}$ and $\xi_{n+p+1}$ are repeated $p+1$ times is called \emph{open} knot vectors.
In what follows, we consider only open knot vectors but the construction is the same for general knot vectors.
Moreover, we may refer a patch as a subdomain and an element as a knot span, i.e. an interval of the form $[\xi_{i},\xi_{i+1}]$.

\begin{figure}[!ht]
\centering
\includegraphics[width=0.6\linewidth]{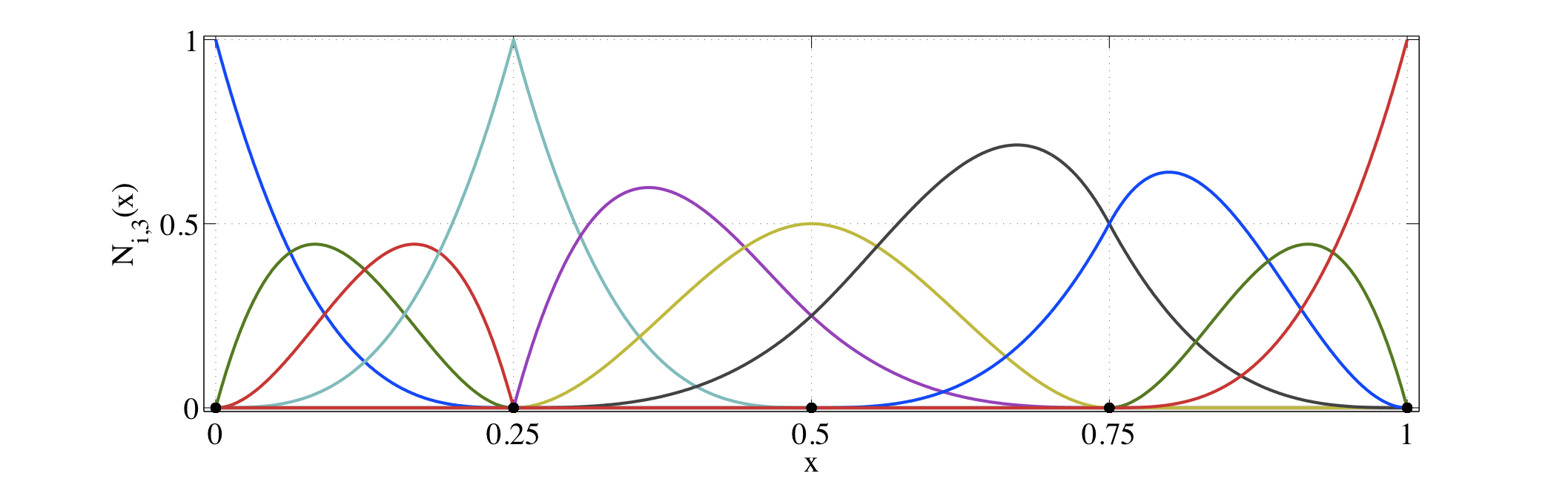}
\caption{Example of B-spline basis functions for $\xi=\left\{0,0,0,0,0.25,0.25,0.25,0.5,0.75,0.75,1,1,1,1\right\}$, $n=10$ and $p=3$. We see that the regularity is related to the multiplicity of each $\xi_{i}$. Moreover, for $\xi_{i}=0.25$ we have $m_{i}=p$ and we see that the basis function is interpolatory at this knot.}
\label{fig:bspline_basis}
\end{figure}

The B-spline functions are constructed recursively with respect to the polynomial order.
For $p=0$ and an open knot vector $\Xi$, we define
\begin{equation}\label{equ:bsplinebasis1}\nonumber
N_{i,0}(x):=\left\{
\begin{array}{ll}
1 & \mbox{if }\xi_{i}\leq x\leq\xi_{i+1},\\
0 & \mbox{otherwise.}
\end{array}
\right.
\end{equation}
For $p=1,2,\dots$, we define recursively the B-spline basis functions as
\begin{equation}\label{equ:bsplinebasis}
N_{i,p}(x):=\frac{x-\xi_{i}}{\xi_{i+p}-\xi_{i}}N_{i,p-1}(x)+\frac{\xi_{i+p+1}-x}{\xi_{i+p+1}-\xi_{i+1}}N_{i+1,p-1}(x).
\end{equation}
We present in Figure \ref{fig:bspline_basis} an example of B-spline basis functions for $n=10$ and $p=3$ and the knot vector $\xi=\left\{0,0,0,0,0.25,0.25,0.25,0.5,0.75,0.75,1,1,1,1\right\}$.
The equation \eqref{equ:bsplinebasis} is called the Cox-de-Boor recursion formula \cite{cox1972numerical,de1972calculating}.
Note that for $p=0,1$, the B-spline basis functions coincide with the FE ones.
The B-splines constitute a partition of the unity, that is
\begin{equation}\label{equ:bsplinepartitionunity}\nonumber
\sum_{i=1}^{n}N_{i,p}(x)=1,\qquad\forall x\in[0,1].
\end{equation}
A second feature is that they are pointwise non-negative, i.e. $N_{i,p}(x)\geq0$, $\forall x\in[0,1]$.
This implies that the coefficients of the mass matrix are greater or equal than zero.
The support of $N_{i,p}$ is $[\xi_{i},\xi_{i+p+1}]$.
The basis function $N_{i,p}$ has $p-m_{i}$ continuous derivatives, where $m_{i}$ is the multiplicity of $\xi_{i}$, i.e. the number of repetitions of $\xi_{i}$.
An important remark is that the B-spline basis functions are not interpolatory at the location of knot values $\xi_{i}$ unless the multiplicity of $\xi_{i}$ is exactly $p$.

We are now in position to define B-spline curves, surfaces and solids in $\mathbb{R}^{d}$.
Let us assume that we are given three sets of B-spline basis functions $\left\{N_{i,p}\right\}$, $\left\{M_{j,q}\right\}$ and $\left\{L_{k,r}\right\}$ constructed on the knot vectors $\left\{\xi_{1},\dots,\xi_{n+p+1}\right\}$, $\left\{\eta_{1},\dots,\eta_{m+q+1}\right\}$ and $\left\{\zeta_{1},\dots,\zeta_{l+r+1}\right\}$ for $i=1,\dots,n$, $j=1,\dots,m$ and $k=1,\dots,l$, respectively.

The B-spline curves are obtained by considering linear combinations of B-spline basis functions.
Let $C_{i}\in\mathbb{R}^{d}$ be the coefficients referred as \emph{control points}, for $i=1,\dots,n$.
We then define a B-spline curve as
\begin{equation}\label{equ:bsplinecurve}\nonumber
S(x):=\sum_{i=1}^{n}N_{i,p}(x)C_{i}.
\end{equation}
Such curves have at least as many continuous derivatives across an element boundary than its underlying B-spline basis function has across the corresponding knot value.
A crucial property of the B-spline curves is that an affine transformation of the curve is obtained by applying the transformation to the control points.
It is the so-called \emph{affine covariance} and play an important role in the affine decomposition \eqref{equ:affinea} when considering parameter-dependent domains.
Now that the univariate B-splines have been introduced, we generalize the definition to higher dimensions by considering a tensor product structure.

Given a so-called \emph{control net} $\left\{C_{i,j}\right\}\subset\mathbb{R}^{d}$ for $i=1,\dots,n$ and $j=1,\dots,m$, we define a B-spline surface as
\begin{equation}\label{equ:bsplinesurf}\nonumber
S(x,y):=\sum_{i=1}^{n}\sum_{j=1}^{m}N_{i,p}(x)M_{j,q}(y)C_{i,j}.
\end{equation}
Several properties of the B-spline surfaces result from their tensor product structures.
For instance, the basis also forms a partition of the unity and the number of continuous partial derivatives are determined from the underlying one-dimensional knot vector and polynomial order.
The local support is also deducted from the one-dimensional basis, that is the support of $N_{i,p}(x)M_{j,q}(y)$ is $[\xi_{i},\xi_{i+p+1}]\times[\eta_{j},\eta_{j+q+1}]$.

Finally, we introduce the definition of a B-spline solid.
Considering a control lattice $\left\{C_{i,j,k}\right\}$ for $i=1,\dots,n$, $j=1,\dots,m$ and $k=1,\dots,l$, it is defined as
\begin{equation}\label{equ:bsplinesolid}\nonumber
S(x,y,z):=\sum_{i=1}^{n}\sum_{j=1}^{m}\sum_{k=1}^{l}N_{i,p}(x)M_{j,q}(y)L_{k,r}(z)C_{i,j,k}.
\end{equation}
The properties of the B-spline solids are a direct extension of those presented in the case of surfaces.
In particular, the affine covariance property still holds for B-spline surfaces and solids.
Note that what has been presented here is valid for a single patch.
The case of multipatches geometries is introduced in the context of NURBS basis functions in the next section.

For the purpose of the analysis presented in Section \ref{sec:NURBS_approx}, we briefly discuss the notions of $h$-refinement, $p$-refinement and $k$-refinement.
A complete discussion can be found in \cite{cottrell2009isogeometric,hughes2005isogeometric}.
The notion of $h$-refinement in FE analysis is similar to the \emph{knot insertion} in IGA.
Let us consider a knot vector $\Xi=\{\xi_{1},\dots,\xi_{n+p+1}\}$ and associated control points $\{B_{1},\dots,B_{n}\}$.
Considering a knot $\bar{\xi}\in[\xi_{k},\xi_{k+1}[$, we then build the new knot vector as $\Xi=\{\xi_{1},\dots,\xi_{k},\bar{\xi},\xi_{k+1},\dots,\xi_{n+p+1}\}$ and the associated control points $\{\bar{B}_{1},\dots,\bar{B}_{n+1}\}$ as
\begin{equation}\label{equ:href_control_points}
\bar{B}_{i}:=\alpha_{i}B_{i}+(1-\alpha_{i})B_{i-1},
\end{equation}
where
\begin{equation}\label{equ:href_alpha}
\alpha_{i}:=\left\{
\begin{array}{ll}
1, & 1\leq i\leq k-p,\\
\frac{\bar{\xi}-\xi_{i}}{\xi_{i+p}-\xi_{i}}, & k-p+1\leq i \leq k,\\
0, & k+1\leq i \leq n+p+2.
\end{array}
\right.
\end{equation}
By choosing the new control points as \eqref{equ:href_control_points} and \eqref{equ:href_alpha}, it is possible to maintain the continuity of the original basis functions.
Note that it is possible to insert repetition of already existing knot values.
This will decrease the regularity of the basis functions at this knot.
An important remark is that the solution spanned by the increased basis functions based contains the one spanned by the original B-splines.
This allows to keep the geometry unchanged by inserting new knots.

The second concept introduced here is the FE $p$-refinement, which analogous is \emph{order elevation}.
It is possible to increase the polynomial order of the basis functions.
To keep the regularity of the previous B-splines, it is necessary to repeat each knot value of the knot vector.
As in the case of knot insertion, the new span contains the one from the original basis functions.

The last notion is the one of $k$-refinement which does not have an analogous in FE analysis.
It is based on the principle that order elevation and knot insertion do not commute.
If we insert a new knot value $\bar{\xi}$, the continuity of the basis functions at this knot will be $C^{p-1}$.
If then we further increase the order of the basis, the multiplicity of $\bar{\xi}$ increase to keep this continuity.
Instead, if we first increase the order of the basis to $q$ and then insert a new knot value, the continuity will be $C^{q-1}$ at this knot.
This second process is called $k$-refinement.
It allows to control the number of new basis functions.
Hence the number of degrees of freedom associated to the B-splines will also be kept under control, which in turn to keep low computational costs.

\subsection{Non-Uniform Rational B-Splines}

The introduction of NURBS allows us to exactly represent domains that it is not possible to describe considering polynomials.
The construction of such geometries in $\mathbb{R}^{d}$ are obtained by projective transformations in $\mathbb{R}^{d+1}$.
It is then possible to construct for instance conic sections.
Such projective transformation yields rational polynomial functions.

The process to construct NURBS basis functions is presented here and follows mainly \cite{cottrell2009isogeometric,hughes2005isogeometric}.
Let us consider a knot vector $\Xi$, a number of basis functions $n$, a polynomial order $p$ and a set of control points $\{B_{i}^{w}\}$ in $\mathbb{R}^{d+1}$ defining a B-spline curve.
Such points are called \emph{projective} control points for the associated NURBS curve.
We then define the control points of the NURBS curve as follows
\begin{eqnarray}
w_{i} & := & (B_{i}^{w})_{d+1},\qquad i=1,\dots,d,\nonumber\\
\left(B_{i}\right)_{j} & := & \left(B_{i}^{w}\right)_{j}/w_{i},\qquad i,j=1,\dots,d,\nonumber
\end{eqnarray}
where $\left(B_{i}\right)_{j}$ is the $j$th component of the vector $B_{i}$.
The scalars $w_{i}$ are called weights.
Let $\{N_{i,p}\}$ be the B-spline basis functions associated to $\Xi$, $n$ and $p$.
Based on the definition of the control points, we can introduce the NURBS basis functions defined as
\begin{equation}\label{equ:nurbs_basis}
R_{i}^{p}(x):=\frac{N_{i,p}(x)w_{i}}{\sum_{i'=1}^{n}N_{i',p}(x)w_{j}}.
\end{equation}
The associated NURBS curve is then defined as
\begin{equation}\label{equ:nurbs_curve}\nonumber
C(x):=\sum_{i=1}^{n}R_{i}^{p}(x)B_{i}.
\end{equation}

Considering the basis functions defined by \eqref{equ:nurbs_basis}, we define NURBS surfaces and solids in the same manner.
To do this, we define rational basis functions for surfaces and solids.
Let $\{M_{j,q}\}$ and $\{L_{k,r}\}$ be B-spline basis functions for $1\leq j\leq m$ and $1\leq k\leq l$.
Moreover, consider projective control nets and lattices $\{B_{i,j}^{w}\}$ and $\{B_{i,j,k}^{w}\}$ in $\mathbb{R}^{d+1}$, respectively.
The weights to construct the NURBS basis functions are given by
\begin{eqnarray}
w_{i,j} & := & \left(B_{i,j}^{w}\right)_{d+1},\qquad i,j=1,\dots,d,\nonumber\\
w_{i,j,k} & := & \left(B_{i,j,k}^{w}\right)_{d+1},\qquad i,j,k=1,\dots,d.\nonumber
\end{eqnarray}
We then define NURBS basis functions for surfaces and solids as
\begin{eqnarray}
R_{i,j}^{p,q}(x,y) & := & \frac{N_{i,p}(x)M_{j,q}(y)w_{i,j}}{\sum_{i'=1}^{n}\sum_{j'=1}^{m}N_{i',p}(x)M_{j',q}(y)w_{i',j'}},\qquad i,j=1\dots,d,\nonumber\\
R_{i,j,k}^{p,q,r}(x,y,z) & := & \frac{N_{i,p}(x)M_{j,q}(y)L_{k,r}(z)w_{i,j,k}}{\sum_{i'=1}^{n}\sum_{j'=1}^{m}\sum_{k'=1}^{l}N_{i',p}(x)M_{j',q}(y)L_{k',r}(z)w_{i',j',k'}},\qquad i,j,k=1\dots,d.\nonumber
\end{eqnarray}

\begin{figure}[!ht]
\centering
\subfloat[][]{\includegraphics[clip=true,width=0.2\linewidth]{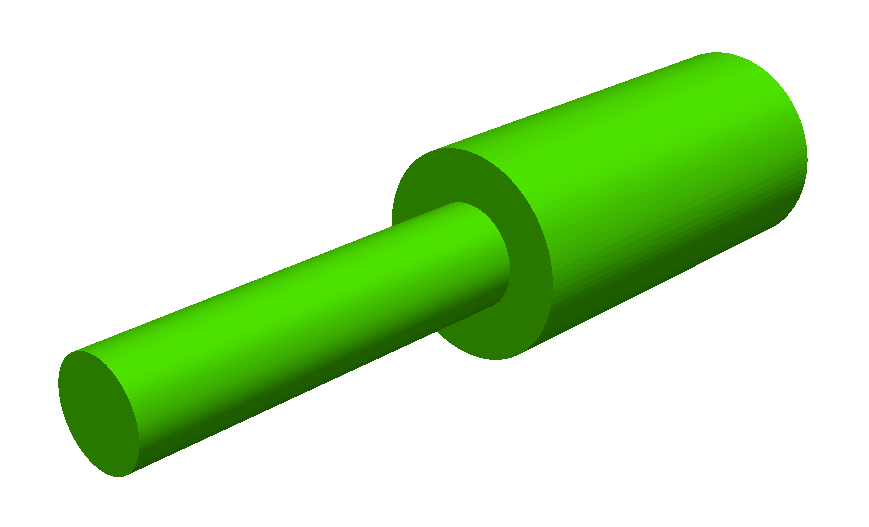}\label{subfig:adhesive_roll}}\
\subfloat[][]{\includegraphics[clip=true,width=0.14\linewidth]{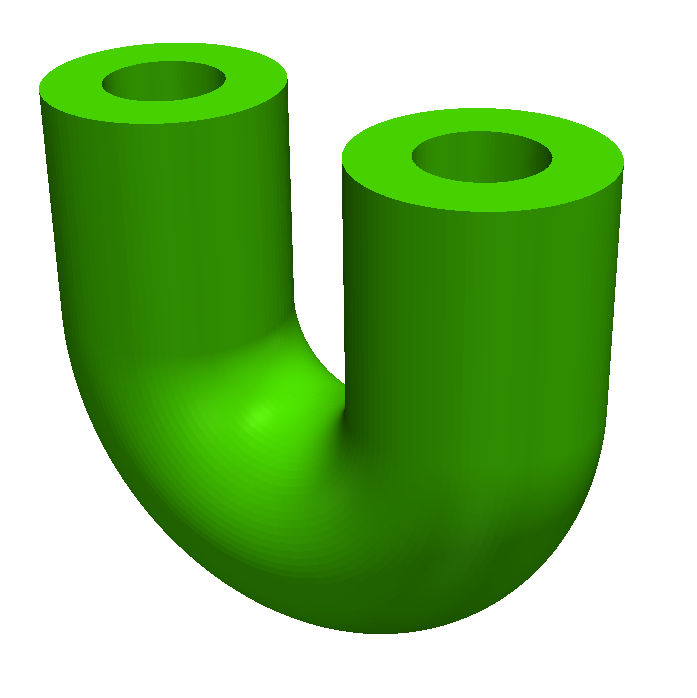}\label{subfig:u_pipe}}\
\subfloat[][]{\includegraphics[clip=true,width=0.2\linewidth]{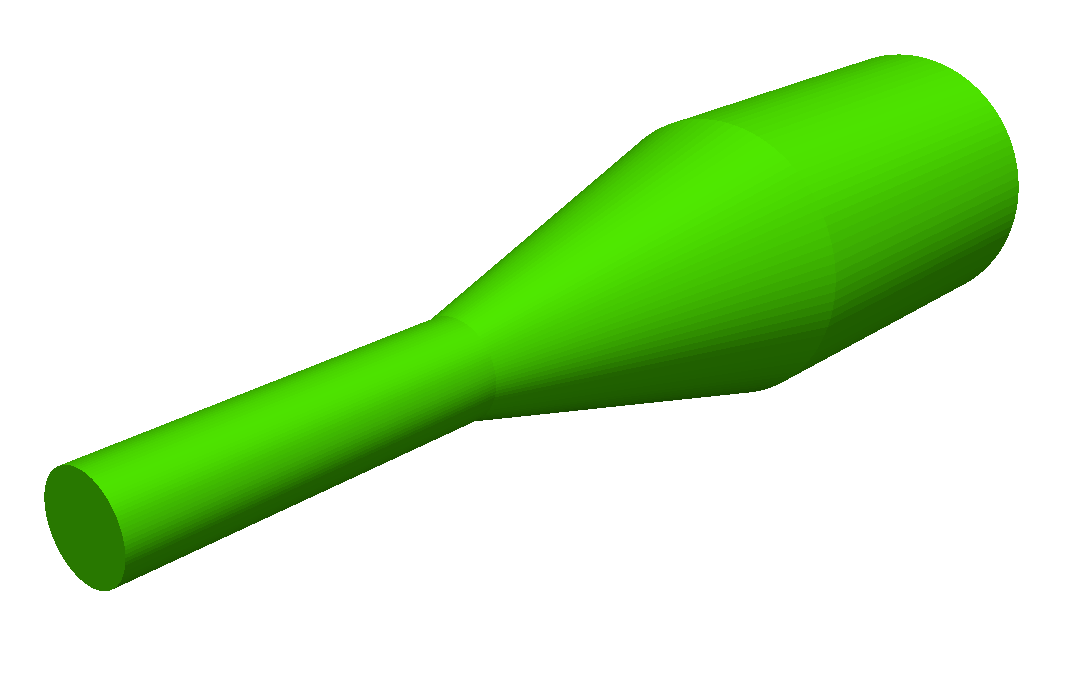}\label{subfig:cylinder_diff_diam}}\
\subfloat[][]{\includegraphics[clip=true,width=0.2\linewidth]{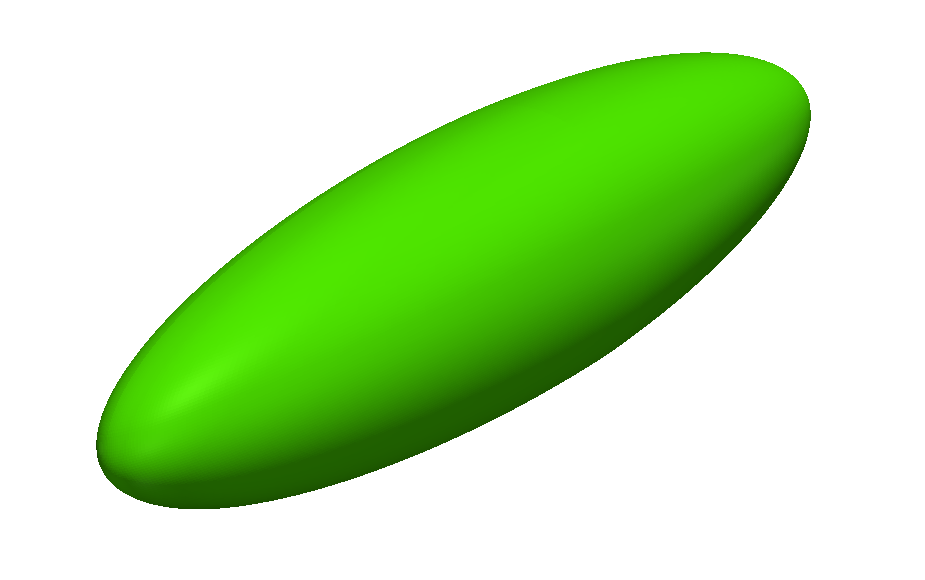}\label{subfig:ellipsoid}}\
\subfloat[][]{\includegraphics[clip=true,width=0.2\linewidth]{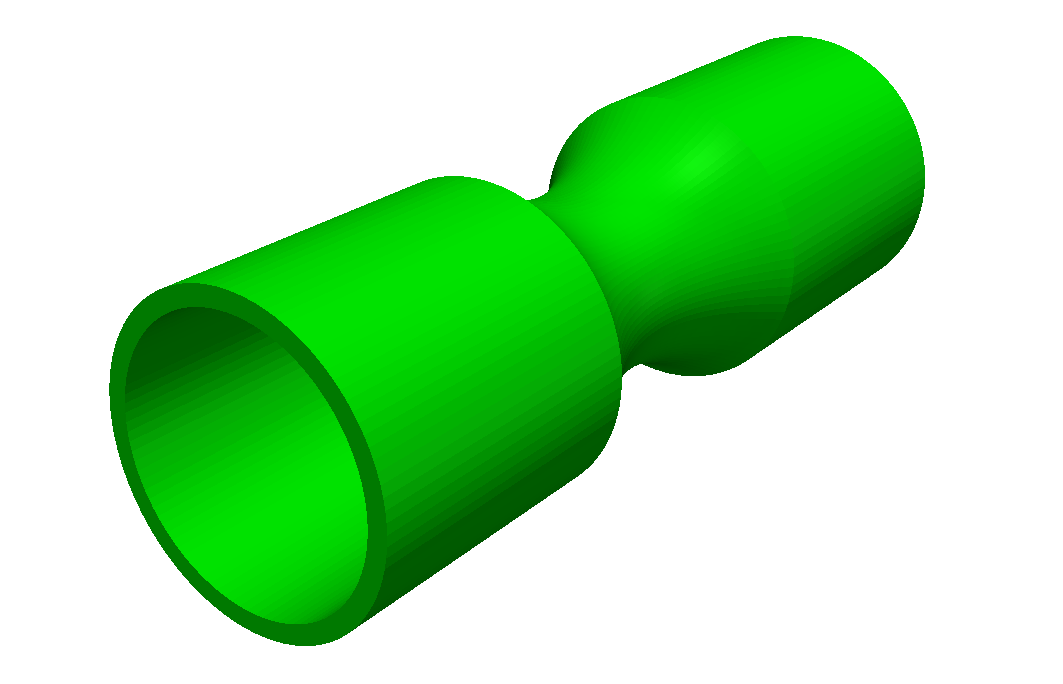}\label{subfig:bended_cylinder}}
\caption{Examples of NURBS solids obtained considering multiple patches. The number of patches used for each example are $3$ \protect\subref{subfig:adhesive_roll}, $3$ \protect\subref{subfig:bended_cylinder}, $3$ \protect\subref{subfig:cylinder_diff_diam}, $4$ \protect\subref{subfig:ellipsoid} and $4$ \protect\subref{subfig:u_pipe}, respectively.}
\label{fig:nurbs_examples}
\end{figure}

The properties stated for the B-spline basis functions also hold for the NURBS.
In particular, they form a partition of the unity and their continuity and support are the same as the underlying B-splines.
The affine covariance property also holds for NURBS functions.
Moreover, the basis functions are interpolatory at knot values where the multiplicity is equal to the order.
The notions of $h$-, $p$- and $k$-refinement generalize to NURBS functions.
Note that if all the weights are equal, the NURBS coincide with the underlying B-splines due to the partition of the unity property.
In nearly all the practical applications, it is necessary to have multiple patches to describe the domain with NURBS functions.
This also allows to have different material properties, each associated to a different patch.
The only feature to pay attention to is the regularity of the basis across the patches interfaces.
Usually, $C^{0}$ is the only regularity guaranteed, but techniques can be used to increase it \cite{cottrell2009isogeometric}.
In Figure \ref{fig:nurbs_examples}, we present several examples of NURBS solids obtained considering multipatches representations.

To simplify the notations, we denote by $R_{i,p}$, $1\leq i\leq n$ the NURBS basis functions and $\{B_{i}\}$ the associated control points for curves, surfaces and solids.
We also use the notation $\Xi$ for the associated knot vectors.
Note that it is a slight abuse of notation because in the case of surfaces and solids, $p$ and $\Xi$ are vectors and matrices, respectively.

For the purpose of our analysis, we require that the computational domain $\Omega$ can be obtained through a NURBS parametrization.
To introduce the notations, we impose that $\Omega$ is parameter-independent.
The parameter-dependent case is treated in the next section.
Let us consider the following decomposition of the domain
\begin{equation}\label{equ:domain_decomp}
\overline{\Omega} = \bigcup_{k=1}^{P_{\text{dom}}}\overline{\Omega}^{k},
\end{equation}
where $\Omega^{j}\cap\Omega^{k}=\emptyset$ for $1\leq j<k\leq P_{\text{dom}}$.
We require that for every subdomain $\Omega^{k}$, there exist $p_{k}$, $n_{k}$, $\Xi_{k}$, NURBS basis functions $\{R_{i,p}^{k}\}$ and associated control points $\mathcal{B}^{k}:=\{B_{i}^{k}\}$ such that for every $y\in\Omega^{k}$, there exists $x\in\mathcal{H}^{d}$ satisfying
\begin{equation}\label{equ:nurbs_representation}
y=F^{k}(x):=\sum_{i=1}^{n_{k}}R_{i,p_{k}}^{k}(x)B_{i}^{k}.
\end{equation}
Here $\mathcal{H}^{d}=[0,1]^{d}$ denotes the unit hypercube in $d$-dimension and $F^{k}:(0,1)^{d}\rightarrow\Omega^{k}$.
Considering for every $1\leq k\leq P_{\text{dom}}$ the function $F^{k}$ defined above, we construct a global mapping $F:(0,1)^{d}\rightarrow\Omega$ which describes the whole computational domain.
We assume that $F$ is smooth and invertible.
In that case, we say that $\Omega$ admits a NURBS representation through $F$.
So far, we have only considered parameter-independent geometries.
In the next section, we introduce the affine preconditionning conditions for parameter-dependent domains.

\subsection{Affine Preconditionning for Parameter-Dependent Domains}\label{subsec:affine_precond}

In many applications, it is of great interest to consider parameter-dependent geometries.
We introduce here the conditions that need to be fulfilled in that case to be able to perform the RB method presented in this paper.
In particular, it is important that the affine decomposition \eqref{equ:affinea} of the bilinear form $a$ still holds.
Let us consider the domain splitting introduced in \eqref{equ:domain_decomp}.
The computational domain for an input parameter $\boldsymbol\mu\in\mathcal{D}$ is denoted $\Omega_{o}(\boldsymbol\mu)$.
Here, the subscript $o$ stands for the original domain.

The domain $\Omega_{o}(\boldsymbol\mu)$ needs to be represented as the image of a reference domain through an affine mapping.
Let us choose $\boldsymbol\mu_{\text{ref}}\in\mathcal{D}$ as a parameter that represents our reference domain, i.e. $\Omega=\Omega_{o}(\boldsymbol\mu_{\text{ref}})$.
Moreover, we denote $\Omega^{k}=\Omega^{k}_{o}(\boldsymbol\mu_{\text{ref}})$ while considering the decomposition \eqref{equ:domain_decomp}.
We need that for every $1\leq k\leq P_{\text{dom}}$, there exists an affine mapping $\mathcal{T}^{k}(\cdot;\boldsymbol\mu):\Omega^{k}\rightarrow\Omega_{o}^{k}(\boldsymbol\mu)$ such that
\begin{equation}\nonumber
\overline{\Omega}_{o}^{k}(\boldsymbol\mu)=\mathcal{T}^{k}(\overline{\Omega}^{k};\boldsymbol\mu).
\end{equation}
The mappings $\mathcal{T}^{k}(\cdot;\boldsymbol\mu)$ have to be bijective and collectively continuous, that is
\begin{equation}\label{equ:affine_map_cont_bij_first}
\mathcal{T}^{k}(x;\boldsymbol\mu)=\mathcal{T}^{l}(x;\boldsymbol\mu),\qquad\forall x\in\overline{\Omega}^{k}\cap\overline{\Omega}^{l},\ 1\leq k<l\leq P_{\text{dom}}.
\end{equation}

Due to the affine covariance property of the NURBS functions, we only need to require that the control points can be obtained as the image of reference control points through an affine mapping.
More formally, let us denote by $\{B_{i}^{k}(\boldsymbol\mu)\}$ the control points associated with the subdomains $\Omega_{o}^{k}(\boldsymbol\mu)$.
We then require that for every $1\leq k\leq P_{\text{dom}}$, there exists an affine mapping $\mathcal{T}^{k}(\cdot;\boldsymbol\mu):\Omega^{k}\rightarrow\Omega_{o}^{k}(\boldsymbol\mu)$ such that
\begin{equation}\nonumber
B_{i}^{k}(\boldsymbol\mu)=\mathcal{T}^{k}(B_{i}^{k};\boldsymbol\mu),
\end{equation}
where $\{B_{i}^{k}\}$ are the control points associated to the reference subdomains $\Omega^{k}$.
Turning to the condition \eqref{equ:affine_map_cont_bij_first}, we require that
\begin{equation}\label{equ:affine_map_cont_bij_second}\nonumber
\mathcal{T}^{k}(B_{i}^{k};\boldsymbol\mu)=\mathcal{T}^{l}(B_{i}^{k};\boldsymbol\mu),\qquad\forall B_{i}^{k}\in\mathcal{B}^{k}\cap\mathcal{B}^{l},\ 1\leq k<l\leq P_{\text{dom}}.
\end{equation}
In other words, we only need to ensure continuity of the mappings through the control points defining the interfaces of patches to obtain the continuity on the whole interface.
More explicitly, we define the affine mappings $\mathcal{T}^{k}$ for every $x\in\overline{\Omega}^{k}$ and $\boldsymbol\mu\in\mathcal{D}$ as
\begin{equation}\nonumber
\mathcal{T}^{k}(x;\boldsymbol\mu):=C^{k}(\boldsymbol\mu)+G^{k}(\boldsymbol\mu)x,
\end{equation}
where $C^{k}:\mathcal{D}\rightarrow\mathbb{R}^{d}$ and $G^{k}:\mathcal{D}\rightarrow\mathbb{R}^{d\times d}$ for every $1\leq k \leq P_{\text{dom}}$.
To define the affine decomposition of the bilinear form $a$, we need to define the Jacobians and inverse of the transformations as
\begin{eqnarray}
J^{k}(\boldsymbol\mu) & := & |\det(G^{k}(\boldsymbol\mu))|\label{equ:defjk},\\
D^{k}(\boldsymbol\mu) & := & \left(G^{k}(\boldsymbol\mu)\right)^{-1},\label{equ:defgk}
\end{eqnarray}
for $1\leq k \leq P_{\text{dom}}$.
Based on the $\mathcal{T}^{k}$ transformations, we can define a global affine mapping $\mathcal{T}:\Omega\rightarrow\Omega_{o}(\boldsymbol\mu)$ as
\begin{equation}\label{equ:global_mapping}\nonumber
\mathcal{T}(x;\boldsymbol\mu):=\mathcal{T}^{k}(x;\boldsymbol\mu),\qquad k=\min\left\{1\leq l \leq P_{\text{dom}}\middle|\ x\in\overline{\Omega}^{l}\right\}.
\end{equation}
The mapping $\mathcal{T}$ is globally bijective and piecewise affine.
The choice of the minimum is arbitrary and could be chosen differently.

\begin{figure}[!ht]
\centering
\subfloat[]{\label{subfig:torus}\includegraphics[clip=true,width=0.24\linewidth]{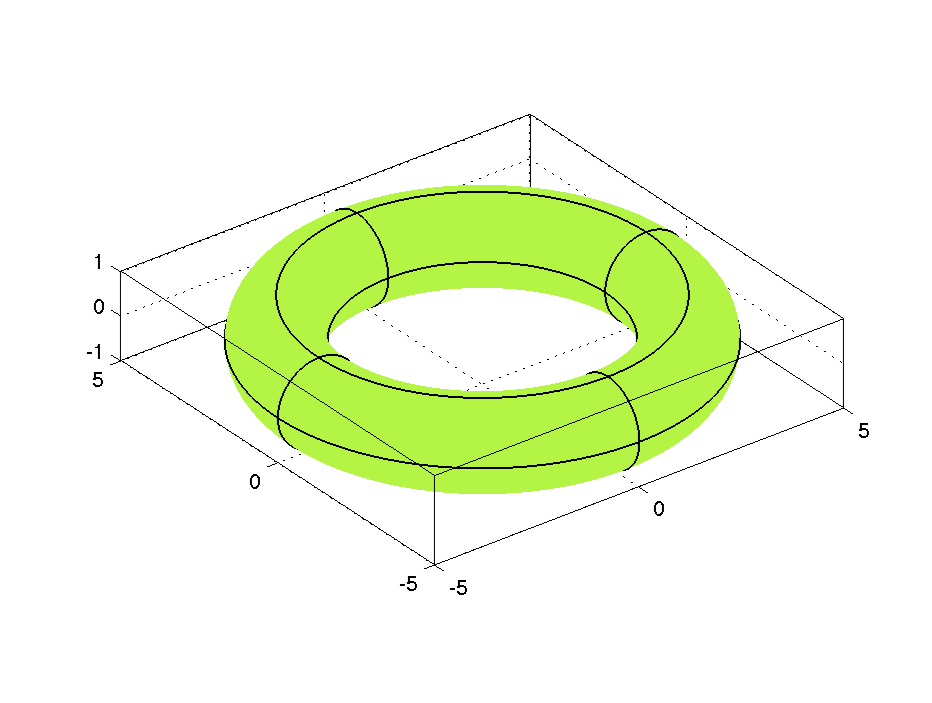}}\
\subfloat[]{\label{subfig:torus_lattice}\includegraphics[clip=true,width=0.24\linewidth]{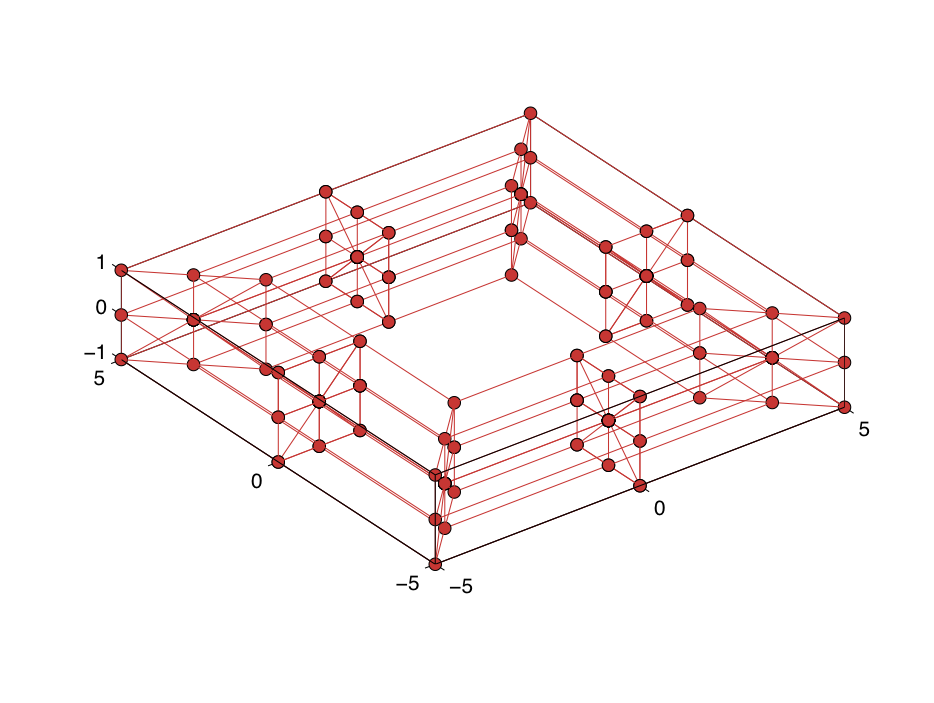}}\
\subfloat[]{\label{subfig:torus_ellipse}\includegraphics[clip=true,width=0.24\linewidth]{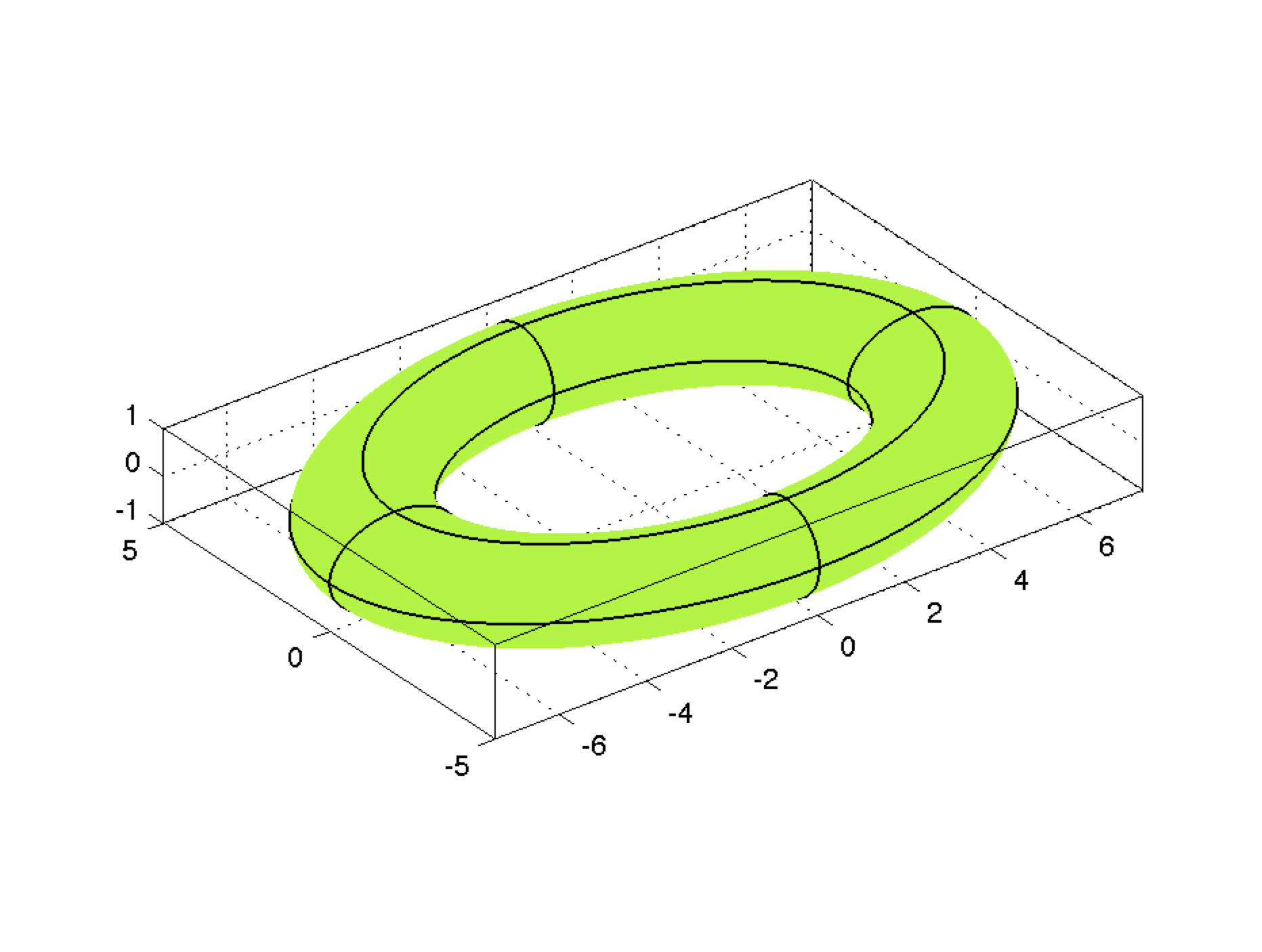}}\
\subfloat[]{\label{subfig:torus_ellipse_lattice}\includegraphics[clip=true,width=0.24\linewidth]{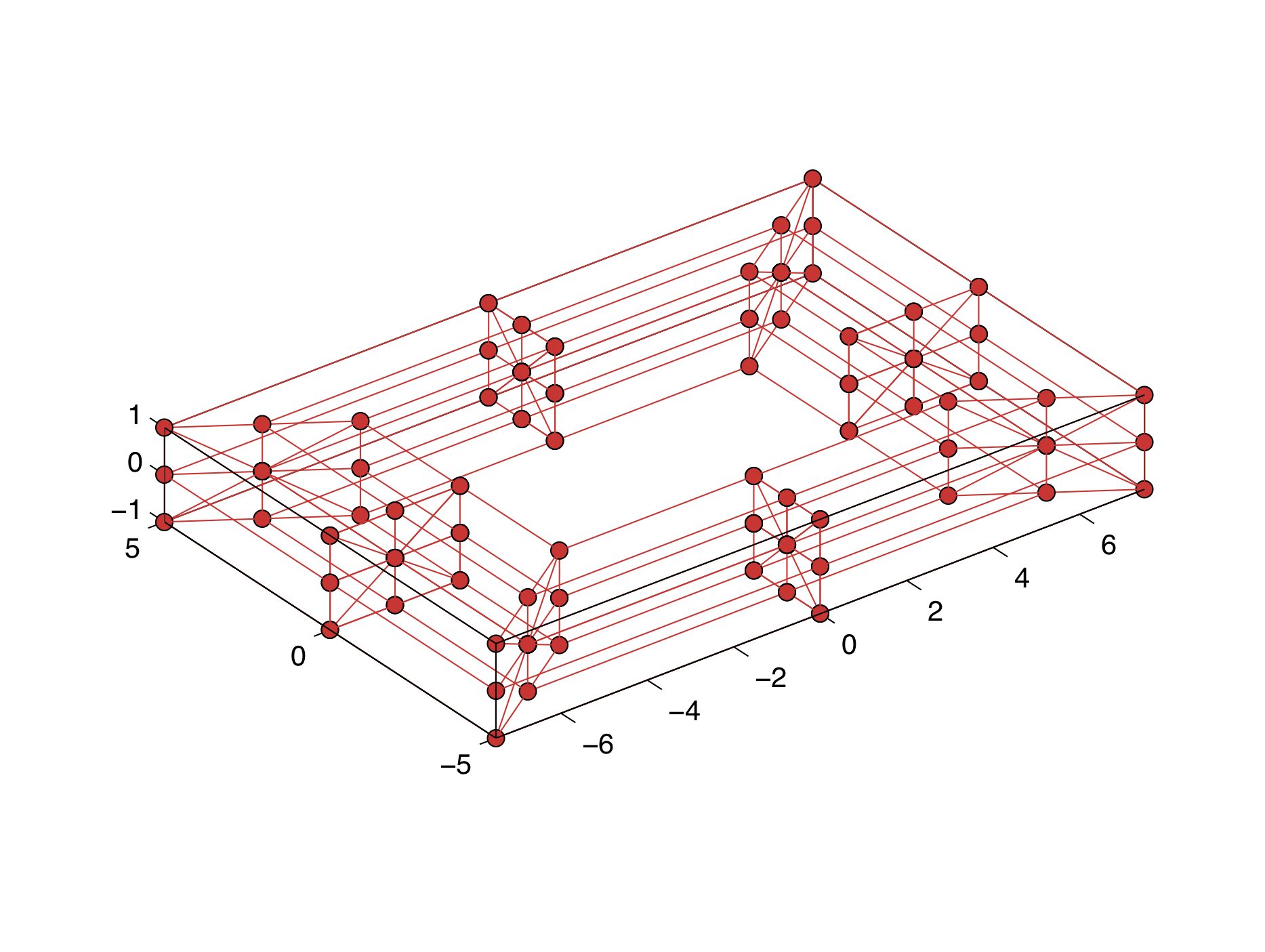}}
\caption{Example of the affine transformation \eqref{equ:transformation_torus} applied to a torus. The original domain and its lattice of control points are presented in \protect\subref{subfig:torus} and \protect\subref{subfig:torus_lattice}, respectively. The torus after transformation is depicted in \protect\subref{subfig:torus_ellipse} while its lattice is presented in \protect\subref{subfig:torus_ellipse_lattice}.}
\label{fig:torus_transformation}
\end{figure}

In what follows, we give an example of an affine transformation applied to a $3$-dimensional toroidal solid.
It is built on four patches, which yields $P_{\text{dom}}=4$, and to every patch are associated $27$ control points.
Based on that, it is possible to uniquely determine $C^{k}$ and $G^{k}$ for $1\leq k\leq P_{\text{dom}}$.
In that case, the transformations are given by
\begin{equation}\label{equ:transformation_torus}
C^{k}=\left(
\begin{array}{cc}
0\\
0\\
0
\end{array}
\right),\qquad
G^{k}=\left(
\begin{array}{ccc}
\mu & 0 & 0\\
0 & 1 & 0\\
0 & 0 & 1
\end{array}
\right),\qquad1\leq k\leq4,
\end{equation}
where we have considered the single parameter $\mu$ that controls the semi-axis $x$.
In Figure \ref{fig:torus_transformation}, the original domain and the transformed one for $\mu=1.5$ are depicted together with their lattices of control points.
We see that the transformation is exactly applied to the control points.

Now that our affine preconditionning assumption has been stated, we need to express our bilinear form on the reference domain.
The problem \eqref{equ:sex}-\eqref{equ:aex} is defined on the original domain $\Omega_{o}(\boldsymbol\mu)$.
To be able to obtain the affine expansion \eqref{equ:affinea} for the bilinear form $a$ arising from the weak formulation of a second-order PDE, we need that the underlying integrals are defined on the reference domain.
This is presented in details in what follows for two-dimensional problems but the case $d=3$ is treated analogously.
On the original domain, the problem is the following one: given a parameter $\boldsymbol\mu\in\mathcal{D}$, evaluate
\begin{equation}\label{equ:sex_original}\nonumber
s_{o}(\boldsymbol\mu)=l(u_{o}(\boldsymbol\mu)),
\end{equation}
where $u_{o}(\boldsymbol\mu)\in X_{o}$ is the solution of
\begin{equation}\label{equ:aex_original}\nonumber
a_{o}(u_{o}(\boldsymbol\mu),v;\boldsymbol\mu)=f_{o}(v),\qquad\forall v\in X.
\end{equation}
Since we are considering second-order partial differential, we require that $a_{o}$ can be written as
\begin{equation}\label{equ:ao_original}\nonumber
a_{o}(v,w;\boldsymbol\mu)=\sum_{k=1}^{P_{\text{dom}}}\int_{\Omega^{k}(\boldsymbol\mu)}\left[
\begin{array}{ccc}
\frac{\partial v}{\partial x_{1}} & \frac{\partial v}{\partial x_{2}} & v
\end{array}
\right]A_{o}^{k}(\boldsymbol\mu)\left[
\begin{array}{c}
\frac{\partial w}{\partial x_{1}}\\
\frac{\partial w}{\partial x_{2}}\\
v
\end{array}
\right],
\end{equation}
where $A^{k}_{o}:\mathcal{D}\rightarrow\mathbb{R}^{3\times 3}$ is a symmetric positive semi-definite matrix.
We express the right-hand side $f_{o}$ in the same way, that is
\begin{equation}\label{equ:fo_original}
f_{o}(v)=\sum_{k=1}^{P_{\text{dom}}}\int_{\Omega^{k}(\boldsymbol\mu)}f_{o}^{k}v,
\end{equation}
where $f_{o}^{k}\in\mathbb{R}$.
Note that it is possible to have a parameter-dependent right-hand side by simply replacing $f_{o}^{k}$ by $f_{o}^{k}(\boldsymbol\mu)$.
As already discussed, we need $\boldsymbol\mu$-independent integrals to be able to fulfill the affine assumption \eqref{equ:affinea} for the bilinear form $a$.
We then consider the problem \eqref{equ:sex}-\eqref{equ:aex} with the bilinear form $a$ expressed as
\begin{equation}\label{equ:a_expansion}
a(v,w;\boldsymbol\mu)=\sum_{k=1}^{P_{\text{dom}}}\int_{\Omega^{k}}\left[
\begin{array}{ccc}
\frac{\partial v}{\partial x_{1}} & \frac{\partial v}{\partial x_{2}} & v
\end{array}
\right]A^{k}(\boldsymbol\mu)\left[
\begin{array}{c}
\frac{\partial w}{\partial x_{1}}\\
\frac{\partial w}{\partial x_{2}}\\
v
\end{array}
\right],
\end{equation}
where the $A^{k}:\mathcal{D}\rightarrow\mathbb{R}^{3\times 3}$ are defined for $\boldsymbol\mu\in\mathcal{D}$ as
\begin{equation}\label{equ:ak_ako}\nonumber
A^{k}(\boldsymbol\mu)=J^{k}(\boldsymbol\mu)\mathcal{G}^{k}(\boldsymbol\mu)A^{k}_{o}(\boldsymbol\mu)\left(\mathcal{G}^{k}(\boldsymbol\mu)\right)^{T},\qquad 1\leq k \leq P_{\text{dom}}.
\end{equation}
Here  the matrices $\mathcal{G}^{k}(\boldsymbol\mu)$ are defined as
\renewcommand\arraystretch{1}
\begin{equation}\label{equ:gk}\nonumber
\mathcal{G}^{k}(\boldsymbol\mu):=\left(
\begin{array}{ccc}
\multicolumn{2}{c}{\multirow{2}{*}{$D^{k}(\boldsymbol\mu)$}} & 0\\
\multicolumn{2}{c}{} & 0 \\
0 & 0 & 1
\end{array}
\right),\qquad 1\leq k \leq P_{\text{dom}},
\end{equation}
where $J^{k}(\boldsymbol\mu)$ and $D^{k}(\boldsymbol\mu)$ are defined by \eqref{equ:defjk} and \eqref{equ:defgk}, respectively.
Note that this holds under the assumptions presented at the begin of this section.
In the same manner, the right-hand side is expressed as
\begin{equation}\label{equ:f_expansion}\nonumber
f(v)=\sum_{k=1}^{P_{\text{dom}}}\int_{\Omega^{k}}f^{k}(\boldsymbol\mu)v,
\end{equation}
where $f^{k}:\mathcal{D}\rightarrow\mathbb{R}$ is defined by
\begin{equation}\label{equ:fk}\nonumber
f^{k}(\boldsymbol\mu)=J^{k}(\boldsymbol\mu)f_{o}^{k},\qquad1\leq k\leq P_{\text{dom}}.
\end{equation}
We can then explicitly expand \eqref{equ:a_expansion} to obtain the affine decomposition \eqref{equ:affinea} for the bilinear form $a$.
In the development presented here, the $A^{k}(\boldsymbol\mu)$ and $f^{k}(\boldsymbol\mu)$ are local to patches and may represent different material properties and geometry variations.

\subsection{Isogeometric Analysis NURBS Approximation of Elliptic Coercive Parametrized PDEs}\label{sec:NURBS_approx}

We present in this section the isogeometrical analysis NURBS approximation of the problem \eqref{equ:sex}-\eqref{equ:aex}.
In this context, the isoparametric concept is considered, that is the solution is represented in the same space as the geometry.
In that case, the mesh of the NURBS is defined as the product of the knot vectors and the elements are the knot spans.
The degrees of freedom associated with the basis functions are called \emph{control variables}.

Let us assume that $\Omega$ admits a NURBS parametrization through $F$ as defined in \eqref{equ:nurbs_representation}.
To simplify the notations, we consider a single set of indices $\{1,\dots,\mathcal{N}\}$ for the degrees of freedom and we write
\begin{equation}\label{equ:nurbs_representation_analysis}
F(x)=\sum_{i=1}^{\mathcal{N}}\tilde{R}_{i,p}(x)B_{i},\qquad x\in(0,1)^{d},
\end{equation}
for NURBS basis functions $\{\tilde{R}_{i,p}\}$ and associated control points $\{B_{i}\}$.
To represent our solution in a finite dimensional space, we need to define the basis functions
\begin{equation}\label{equ:nurbs_basis_comp_space}
R_{i,p}:=\tilde{R}_{i,p}\circ F^{-1},
\end{equation}
where $F$ is the invertible mapping defined by \eqref{equ:nurbs_representation_analysis}.
Based on that representation, we construct the NURBS approximation space
\begin{equation}\label{equ:nurbs_approx_space}
X^{\mathcal{N}}:=\text{span}\left\{R_{i,p}\right\}_{1\leq i\leq\mathcal{N}}\subset X.
\end{equation}
As already discussed in Section \ref{subsec:bsplines}, the process of knot insertion does not change the underlying geometry.
In that setting, increasing $\mathcal{N}$ does not change the shape of the parametrized domain and so we keep the exact parametrization while refining the mesh.
For approximation properties of NURBS approximation spaces, we refer the reader to \cite{bazilevs2006isogeometric}.
We approximate the solution of \eqref{equ:sex}-\eqref{equ:aex} by an element of $X^{\mathcal{N}}$.
The approximate problem is the following one: given a parameter $\boldsymbol\mu\in\mathcal{D}$, evaluate
\begin{equation}\label{equ:s_nurbs}
s^{\mathcal{N}}(\boldsymbol\mu)=l(u^{\mathcal{N}}(\boldsymbol\mu)),
\end{equation}
where $u^{\mathcal{N}}(\boldsymbol\mu)\in X^{\mathcal{N}}$ is the solution of
\begin{equation}\label{equ:a_nurbs}
a(u^{\mathcal{N}}(\boldsymbol\mu),v;\boldsymbol\mu)=f(v),\qquad\forall v\in X^{\mathcal{N}}.
\end{equation}
Considering the basis $\{R_{i,p}\}$ for $X^{\mathcal{N}}$ defined by \eqref{equ:nurbs_basis_comp_space}, we extend the NURBS solution $u^{\mathcal{N}}(\boldsymbol\mu)$ for $\boldsymbol\mu\in\mathcal{D}$ as
\begin{equation}\label{equ:nurbs_sol_expansion}\nonumber
u^{\mathcal{N}}(x,\boldsymbol\mu)=\sum_{i=1}^{\mathcal{N}}u^{\mathcal{N}}_{i}(\boldsymbol\mu)R_{i,p}(x),\qquad x\in\Omega,
\end{equation}
where the coefficients $u^{\mathcal{N}}_{i}(\boldsymbol\mu)$ are called control variables.
The regularity of $u^{\mathcal{N}}(\boldsymbol\mu)$ follows from that of the NURBS basis.
For instance, the continuity of the solution across element boundaries depends on the continuity of the underlying basis functions across the associated knot span.

Our goal then becomes to solve the problem \eqref{equ:s_nurbs}-\eqref{equ:a_nurbs} with high precision.
However, for real-time context and many query problems, it would be computationally unaffordable to approximate the solution for each input parameter.
For that reason, we introduce in the next section a method to approximate such solution with reduced computation costs.

\section{Reduced Basis Method for Isogeometric Analysis NURBS Approximation}\label{sec:rb}

As it has already been pointed out, it is computationally unaffordable to compute a new NURBS solution for every input parameter $\boldsymbol\mu$.
The goal of the RB method is then to approximate the NURBS solution $u^{\mathcal{N}}(\boldsymbol\mu)$ with reduced computational costs.
Considering $\mathcal{N}$ sufficiently large, we have that $u^{\mathcal{N}}(\boldsymbol\mu)$ is close enough to $u(\boldsymbol\mu)$ in a certain norm so that the NURBS approximation can be viewed as the "truth" solution.

Given a positive integer $N_{\text{max}}\ll\mathcal{N}$, we construct a sequence of approximation spaces
\begin{equation}\label{equ:nestedrbspaces}
X_{1}^{\mathcal{N}}\subset X_{2}^{\mathcal{N}}\subset\dots\subset X_{N_{\text{max}}}^{\mathcal{N}}\subset X^{\mathcal{N}}.
\end{equation}
Those spaces are obtained considering a Greedy algorithm presented more in details in Section \ref{subsec:greedy}.
The hierarchical hypothesis \eqref{equ:nestedrbspaces} is important to ensure the efficiency of the method.
Several spaces can be considered to construct such sequence, but they all focus on the smooth parametric manifold $\mathcal{M}^{\mathcal{N}}:=\left\{u^{\mathcal{N}}(\boldsymbol\mu)\middle |\ \boldsymbol\mu\in\mathcal{D}\right\}$.
If it is smooth enough, we can expect it to be well approximated by low-dimensional spaces.
In what follows, we consider the special case of Lagrange reduced basis spaces built using a master set of parameter points $\boldsymbol\mu^{n}\in\mathcal{D}$, $1\leq n\leq N_{\text{max}}$.
Other examples such as the POD spaces \cite{rozza2008reduced} could be considered.
For $1\leq N\leq N_{\text{max}}$, we define $S^{N}:=\left\{\boldsymbol\mu^{1},\dots,\boldsymbol\mu^{N}\right\}$ and the associated Lagrange RB spaces
\begin{equation}\label{equ:lagspaces}\nonumber
X_{N}^{\mathcal{N}}:=\text{span}\left\{u^{\mathcal{N}}(\boldsymbol\mu^{n})\middle|\ 1\leq n\leq N\right\}.
\end{equation}
The selection of the snapshots $u^{\mathcal{N}}(\boldsymbol\mu^{n})$ is one of the crucial points of the RB method and is further investigated in the next section.
We apply the Gram-Schmidt process in the $(\cdot,\cdot)_{X}$ inner product to the snapshots $u^{\mathcal{N}}(\boldsymbol\mu^{n})$ in order to obtain mutually orthonormal functions $\zeta^{\mathcal{N}}_{n}$.
In that case, we have $X_{N}^{\mathcal{N}}=\text{span}\left\{\zeta^{\mathcal{N}}_{n}\middle|\ 1\leq n\leq N\right\}$.
Since colinearities are avoided using the Gram-Schmidt process, we are ensured that the $N$ obtained is minimal.
The RB approximation of the problem \eqref{equ:s_nurbs}-\eqref{equ:a_nurbs} is obtained considering Galerkin projection: given $\boldsymbol\mu\in\mathcal{D}$, evaluate
\begin{equation}\label{equ:srb}
s^{\mathcal{N}}_{N}(\boldsymbol\mu)=f(u^{\mathcal{N}}_{N}(\boldsymbol\mu)),
\end{equation}
where $u^{\mathcal{N}}_{N}(\boldsymbol\mu)\in X^{\mathcal{N}}_{N}$ is the solution of
\begin{equation}\label{equ:weakenergydensityrb}
a(u^{\mathcal{N}}_{N}(\boldsymbol\mu),v;\boldsymbol\mu)=f(v),\qquad\forall v\in X^{\mathcal{N}}_{N}.
\end{equation}
Since the particular compliant case is considered, we obtain
\begin{equation}\label{equ:outputerror}
s^{\mathcal{N}}(\boldsymbol\mu)-s^{\mathcal{N}}_{N}(\boldsymbol\mu)=\|u^{\mathcal{N}}(\boldsymbol\mu)-u^{\mathcal{N}}_{N}(\boldsymbol\mu)\|_{\boldsymbol\mu}^{2},
\end{equation}
where $\|\cdot\|_{\boldsymbol\mu}$ is the energy norm induced by the inner product $a(\cdot,\cdot;\boldsymbol\mu)$.
In Section \ref{subsec:aposteriori}, we present an example of an inexpensive and efficient a posteriori error estimator $\Delta_{N}(\boldsymbol\mu)$ for $\|u^{\mathcal{N}}(\boldsymbol\mu)-u^{\mathcal{N}}_{N}(\boldsymbol\mu)\|_{\boldsymbol\mu}$ on which the Greedy algorithm is based.
Due to the relation \eqref{equ:outputerror}, it is possible to ensure that the error arising from the RB approximation on our output of interest is bounded by a prescribed tolerance.

Since $u^{\mathcal{N}}_{N}(\boldsymbol\mu)\in X^{\mathcal{N}}_{N}=\text{span}\left\{\zeta^{\mathcal{N}}_{n}\middle|\ 1\leq n\leq N\right\}$, we expand it as
\begin{equation}\label{equ:extensionen}
u^{\mathcal{N}}_{N}(\boldsymbol\mu)=\sum_{m=1}^{N}u^{\mathcal{N}}_{N,m}(\boldsymbol\mu)\zeta^{\mathcal{N}}_{m}.
\end{equation}
The unknowns then become the coefficients $u^{\mathcal{N}}_{N,m}(\boldsymbol\mu)$.
Inserting \eqref{equ:extensionen} in \eqref{equ:srb} and \eqref{equ:weakenergydensityrb} and using the hypothesis that $f$ is linear and $a$ bilinear, we obtain
\begin{equation}\label{equ:extensionsrb}\nonumber
s^{\mathcal{N}}_{N}(\boldsymbol\mu)=\sum_{m=1}^{N}u^{\mathcal{N}}_{N,m}(\boldsymbol\mu)f(\zeta^{\mathcal{N}}_{m}),
\end{equation}
and
\begin{equation}\label{equ:extensionweakenergydensityrb}
\sum_{m=1}^{N}u^{\mathcal{N}}_{N,m}(\boldsymbol\mu)a(\zeta^{\mathcal{N}}_{m},\zeta^{\mathcal{N}}_{n};\boldsymbol\mu)=f(\zeta^{\mathcal{N}}_{n}),\qquad1\leq n\leq N.
\end{equation}
The stiffness matrix associated to the system \eqref{equ:extensionweakenergydensityrb} is of size $N\times N$ with $N\leq N_{\text{max}}\ll\mathcal{N}$.
It yields a considerably smaller computational effort than to solve the system associated to \eqref{equ:a_nurbs}, which matrix is of size $\mathcal{N}\times\mathcal{N}$.
However, the formation of the stiffness matrix involves the computation of the $\zeta^{\mathcal{N}}_{m}$ associated with the $\mathcal{N}$-dimensional NURBS space.

This drawback is avoided by constructing an Offline-Online procedure taking advantage of the affine decomposition \eqref{equ:affinea}.
In the Offline stage, the Greedy algorithm is used to construct the set of parameters $S^{N}$.
Then, the $u^{\mathcal{N}}(\boldsymbol\mu^{n})$ and the $\zeta^{\mathcal{N}}_{n}$ are built for $1\leq n\leq N$.
The $f(\zeta^{\mathcal{N}}_{n})$ and $a^{q}(\zeta^{\mathcal{N}}_{m},\zeta^{\mathcal{N}}_{n})$ are also formed and stored.
Note the importance here of the affine decomposition.
It implies that the vector and matrices stored are independent of the input parameter $\boldsymbol\mu$.

In the Online part, the stiffness matrix associated to \eqref{equ:extensionweakenergydensityrb} is assembled considering the affine decomposition \eqref{equ:affinea}.
This yields
\begin{equation}\nonumber
a(\zeta^{\mathcal{N}}_{m},\zeta^{\mathcal{N}}_{n};\boldsymbol\mu)=\sum_{q=1}^{Q_{a}}\Theta^{q}_{a}(\boldsymbol\mu)a^{q}(\zeta^{\mathcal{N}}_{m},\zeta^{\mathcal{N}}_{n}),\qquad1\leq m,n\leq N.
\end{equation}
The same process is applied to the right-hand side $f$.
The $N\times N$ system \eqref{equ:weakenergydensityrb} is then solved to obtain $u^{\mathcal{N}}_{N,m}(\boldsymbol\mu)$, $1\leq m\leq N$.
Finally, the output of interest \eqref{equ:srb} is computed considering the coefficients obtained.

As already discussed, one of the main feature of the RB method is that we have a posteriori error estimators $\Delta_{N}(\boldsymbol\mu)$ for $\|u^{\mathcal{N}}(\boldsymbol\mu)-u^{\mathcal{N}}_{N}(\boldsymbol\mu)\|_{\boldsymbol\mu}^{2}=s^{\mathcal{N}}(\boldsymbol\mu)-s^{\mathcal{N}}_{N}(\boldsymbol\mu)$ whose computation costs are independent of $\mathcal{N}$.
It allows us to certify our method and make it reliable.
A discussion on such estimators is presented in Section \ref{subsec:aposteriori}.

\subsection{Greedy Algorithm for the Snapshots Selection}\label{subsec:greedy}

One of the most important step taking place in the Offline stage is the selection of the parameters $\boldsymbol\mu^{n}$, $1\leq n\leq N$.
Several algorithms are available in the literature \cite{rozza2008reduced} but we introduce here a greedy procedure for completeness.
The general idea of this procedure is to retain at iteration $N$ the snapshot $u^{\mathcal{N}}(\boldsymbol\mu^{N})$ which approximation by $X^{\mathcal{N}}_{N-1}$ is the worst.
Let us assume that we are given a finite sample of points $\Xi\subset\mathcal{D}$ and pick randomly a first parameter $\boldsymbol\mu^{1}\in\Xi$.
Then for $N=2,\dots,N_{\text{max}}$, compute
\begin{equation}\nonumber
\boldsymbol\mu^{N}:=\arg\max_{\boldsymbol\mu\in\Xi}\Delta_{N-1}(\boldsymbol\mu),
\end{equation}
where $\Delta_{N}(\boldsymbol\mu)$ is a sharp and inexpensive a posteriori error estimator for $\|u^{\mathcal{N}}(\boldsymbol\mu)-u^{\mathcal{N}}_{N}(\boldsymbol\mu)\|_{H^{1}_{0}(\Omega)}$ or $\|u^{\mathcal{N}}(\boldsymbol\mu)-u^{\mathcal{N}}_{N}(\boldsymbol\mu)\|_{\boldsymbol\mu}$.
The algorithm is typically stopped when $\Delta_{N}(\boldsymbol\mu)$ is smaller than a prescribed tolerance for every $\boldsymbol\mu\in\Xi$.
It is clear that the precision of the approximation spaces obtain increase with the size of the sample considered.

Since $X_{N-1}^{\mathcal{N}}\subset X_{N}^{\mathcal{N}}$, we expect to have $\Delta_{N}(\boldsymbol\mu)\leq\Delta_{N-1}(\boldsymbol\mu)$, which ensures that $N_{\text{max}}<\infty$.
Even if this procedure has not been proven to convergence, it is widely used and many examples have been presented to illustrate its convergence.
The derivation of $\Delta_{N}(\boldsymbol\mu)$ is crucial for the Greedy and we introduce an example of such estimator in the next section.

\subsection{A posteriori error estimators for elliptic coercive partial differential equations}\label{subsec:aposteriori}

The main ingredient of the Greedy algorithm procedure is the computation of the error estimator, which has to be independent of $\mathcal{N}$.
In fact, it is used online to certify that the error of our RB approximation with respect to the truth solution is under control.
For completeness, the derivation of such estimator is presented here when the so-called compliant case is considered, i.e. $a$ is symmetric and $f=l$.
See e.g. \cite{rozza2008reduced} for the non-compliant case.
Let us introduce the error $e^{\mathcal{N}}(\boldsymbol\mu)=u^{\mathcal{N}}(\boldsymbol\mu)-u^{\mathcal{N}}_{N}(\boldsymbol\mu)\in X^{\mathcal{N}}$ which satisfies the following equation
\begin{equation}\label{equ:resid}
a(e^{\mathcal{N}}(\boldsymbol\mu),v;\boldsymbol\mu)=f(v;\boldsymbol\mu)-a(u^{\mathcal{N}}_{N}(\boldsymbol\mu),v;\boldsymbol\mu)=:r(v;\boldsymbol\mu),\ \forall v\in X^{\mathcal{N}}
\end{equation}
where $r(\cdot;\boldsymbol\mu)\in\left(X^{\mathcal{N}}\right)^{\prime}$ is the residual and $\left(X^{\mathcal{N}}\right)^{\prime}$ denotes the dual space of $X^{\mathcal{N}}$.
To define our a posteriori error estimator, we need to have a lower bound $\alpha_{\text{LB}}^{\mathcal{N}}(\boldsymbol\mu)$ of $\alpha^{\mathcal{N}}(\boldsymbol\mu)$ such that $0<\alpha_{\text{LB}}^{\mathcal{N}}(\boldsymbol\mu)\leq\alpha^{\mathcal{N}}(\boldsymbol\mu)$ $\forall\boldsymbol\mu\in\mathcal{D}$ and the online costs to compute $\alpha_{\text{LB}}^{\mathcal{N}}(\boldsymbol\mu)$ are independent of $\mathcal{N}$.
We then define the following a posteriori error estimator
\begin{equation}\label{equ:aposteriori}\nonumber
\Delta_{N}(\boldsymbol\mu):=\frac{\|r(\cdot;\boldsymbol\mu)\|_{\left(X^{\mathcal{N}}\right)^{\prime}}}{\alpha_{\text{LB}}^{\mathcal{N}}(\boldsymbol\mu)}.
\end{equation}

To compute $\|r(\cdot;\boldsymbol\mu)\|_{\left(X^{\mathcal{N}}\right)^{\prime}}$, the main ingredients are to use the affine assumption \eqref{equ:affinea} on $a$ and the expansion \eqref{equ:extensionen} of $u_{N}^{\mathcal{N}}(\boldsymbol\mu)$ in the space $X_{N}^{\mathcal{N}}$.
Then, using the definition \eqref{equ:resid} of the residual, this leads to a system depending only on $N$ for every $\boldsymbol\mu$, which makes the computation independent of $\mathcal{N}$.

The procedure used to compute the coercivity lower bound $\alpha_{\text{LB}}^{\mathcal{N}}(\boldsymbol\mu)$ is the so-called successive constraint method (SCM) \cite{scm}.
Considering sets based on parameter samples and the terms $\Theta^{a}_{q}$ of the affine decomposition \eqref{equ:affinea}, it is possible to reduce this problem to a linear optimization problem.
This method works by taking into account neighbour informations for the parameters and its precision increases with the size of the neighbourhood considered.
The SCM also creates a coercivity upper bound $\alpha_{\text{UB}}^{\mathcal{N}}(\boldsymbol\mu)$ of $\alpha^{\mathcal{N}}(\boldsymbol\mu)$ in the same manner.
The algorithm is stopped when $\max_{\boldsymbol\mu\in\Xi}\left(\alpha_{\text{UB}}^{\mathcal{N}}(\boldsymbol\mu)-\alpha_{\text{LB}}^{\mathcal{N}}(\boldsymbol\mu)/\alpha_{\text{UB}}^{\mathcal{N}}(\boldsymbol\mu)\right)$ is smaller than a prescribed tolerance $\varepsilon$.

We emphasize on the fact that it is very important that the costs associated to the computation of $\Delta_{N}$ are independent of $\mathcal{N}$.
That allows us to develop the Offline-Online procedure discussed in Section \ref{sec:rb}, which is a crucial ingredient for the reduced basis method.

\section{Numerical Illustrations}\label{sct:numerics}

In this section, we present several numerical illustrations of the method introduced in this paper.
The first example considered is a case of heat conduction involving only physical parameters, i.e. we use different conductivity coefficients in regions of the domain.
Then a case containing geometrical parameters is introduced.
The aim of the first two illustrations is to present the possibilities that are allowed while using NURBS basis functions.
For this reason, both are computed over curvy three dimensional domains.
In particular, the second example illustrates the theory developed in Section \ref{subsec:affine_precond} to treat parameter dependent geometries.
All computations have been performed using the \texttt{Matlab} \cite{matlab2013} packages \texttt{GeoPDEs} \cite{de2011geopdes} and \texttt{rbMIT} \cite{huynh126002007} for the NURBS and RB approximations, respectively.

Our goal in this section is to present standard examples to show that the method under consideration yields indeed good results.
For this reason, all cases involve simple elliptic equations of the form
\begin{equation}\label{equ:elliptic_numerics}
\begin{array}{rcl}
-\nabla\left(a(\boldsymbol\mu,x)\nabla u\right) & = f, & \mbox{in }\Omega(\boldsymbol\mu),\\
u & = g, & \mbox{on }\Gamma_D(\boldsymbol\mu),\\
a(\boldsymbol\mu,x)\frac{\partial u}{\partial n} & = h, & \mbox{on }\Gamma_{N}(\boldsymbol\mu),
\end{array}
\end{equation}
with $\Gamma_{N}\cap\Gamma_{D}=\emptyset$ and $\partial\Omega=\Gamma_{N}\cup\Gamma_{D}$.
In particular, note that $f$, $g$ and $h$ are parameter independent.
Moreover, we only deal with piecewise constant conductivity coefficients $a(\boldsymbol\mu,x)$.
For all examples, the prescribed tolerance for the greedy algorithm is $10^{-6}$.

\subsection{Physical Parameters for Heat Conduction in a Pipeline}\label{subsec:pipeline}

\begin{figure}[!ht]
\centering
\includegraphics[clip=true,width=0.35\linewidth]{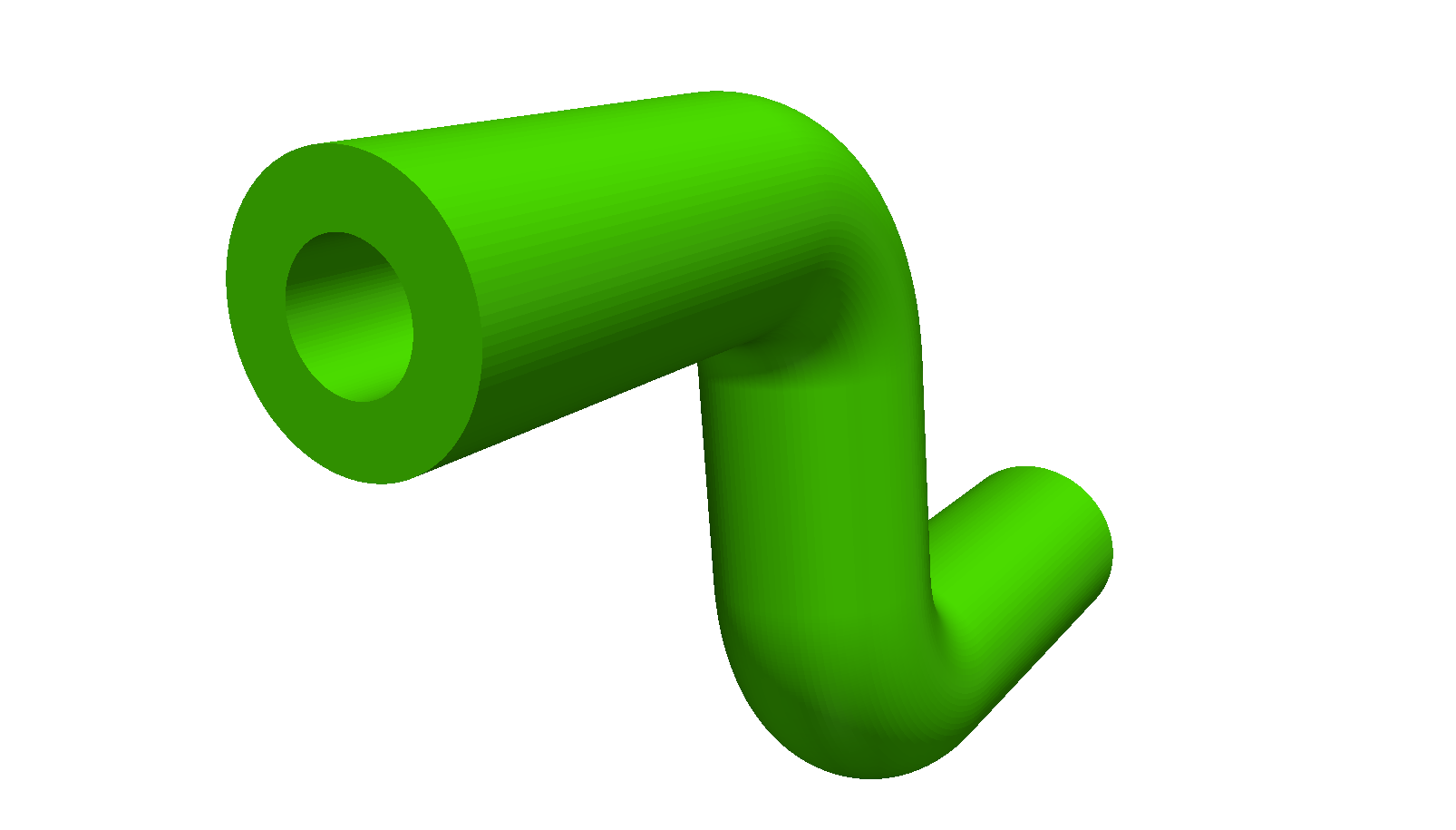}
\caption{Computation domain for the pipeline test case. Five patches were necessary to build the structure, one for each straight part and one for each angle.}
\label{fig:domain_pipeline}
\end{figure}

In this first example, we consider heat conduction in a pipeline.
The domain under consideration is depicted in Figure \ref{fig:domain_pipeline}.
It is built on $5$ different patches, one for every straight part and one for each of the curvy ones.
The domain is parameter independent and we consider three parameters $\boldsymbol\mu=\left(\mu_{1},\mu_{2},\mu_{3}\right)\in[1,5]^{3}$, each one being the conductivity coefficient in one of the straight portion.
More precisely,
\begin{equation*}
a(\boldsymbol\mu,x):=\mu_{1}\chi_{\Omega_{1}}+\chi_{\Omega_{2}}+\mu_{2}\chi_{\Omega_{3}}+\chi_{\Omega_{4}}+\mu_{3}\chi_{\Omega_{5}},
\end{equation*}
where $\chi_{\Omega_{i}}$ is the characteristic function over the $i$th patch, $1\leq i\leq5$.
Let us denote the input boundary by $\Gamma_{\text{in}}$, the output by $\Gamma_{\text{out}}$ and the inner and outer circular ones by $\Gamma_{\text{curve}}$.
The functions are given by $f=0$, $g=0$ and $h=\chi_{\Gamma_{\text{out}}}$ and the associated boundaries are given by $\Gamma_{D}:=\Gamma_{\text{in}}$ and $\Gamma_{N}:=\Gamma_{\text{curve}}\cup\Gamma_{\text{out}}$.
This simulation can be interpreted as heat conduction in a metal pipe where different metals constitute the structure and an imposed temperature is considered on one of the flat faces.

\begin{figure}[!ht]
\centering
\includegraphics[clip=true,width=0.35\linewidth]{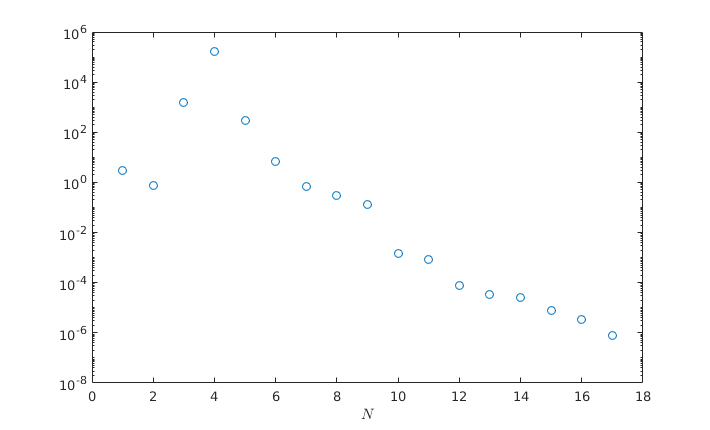}
\caption{Convergence of the greedy algorithm (see Section \ref{subsec:greedy}) for the pipeline test case of Section \ref{subsec:pipeline}.}
\label{fig:conv_greedy_pipeline}
\end{figure}

The number of degrees of freedom for the NURBS approximation is $\mathcal{N}=16650$ while the size of the RB space is $17$, which yields a big reduction of the computational costs.
The computation time to perform the offline step is $27$ minutes and the average evaluation time for the RB approximation is $5\cdot10^{-4}$ seconds.
We present in Figure \ref{fig:conv_greedy_pipeline} the convergence of the greedy algorithm.

\begin{figure}[!ht]
\centering
\subfloat[]{\label{subfig:pipeline_1_1_1}\includegraphics[clip=true,width=0.32\linewidth]{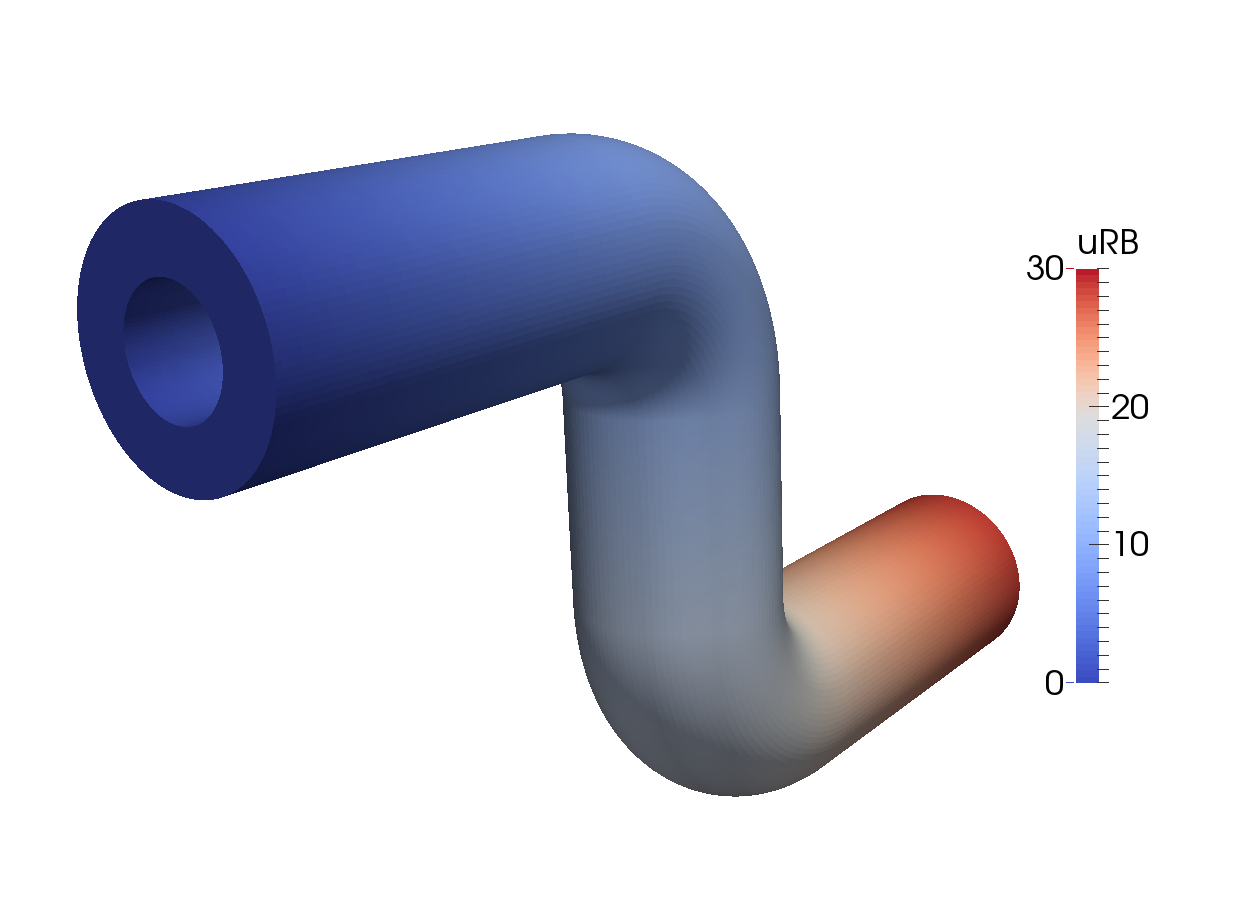}}\
\subfloat[]{\label{subfig:pipeline_3_2_5}\includegraphics[clip=true,width=0.32\linewidth]{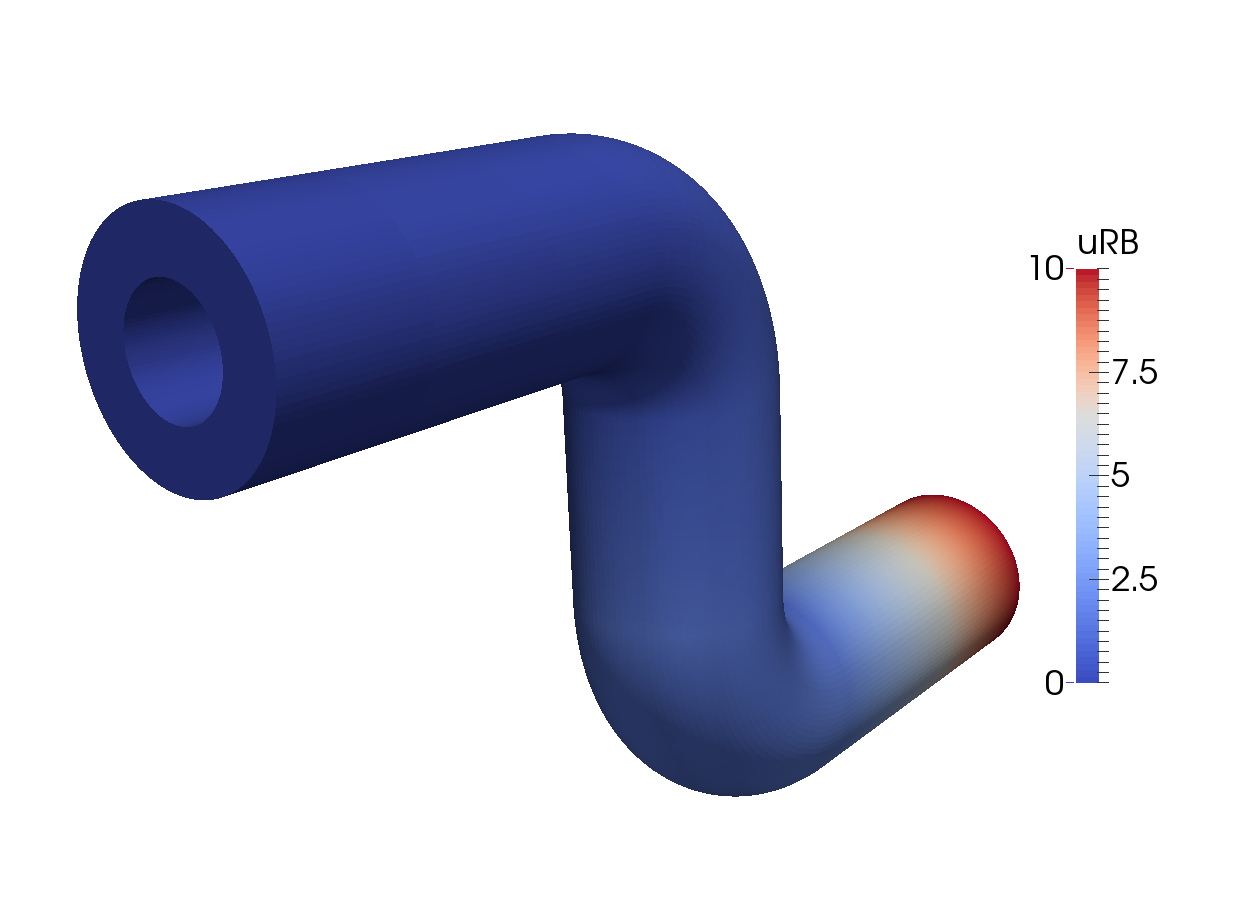}}
\caption{Reduced basis approximation of the pipeline test case for the physical parameters \protect\subref{subfig:pipeline_1_1_1} $\boldsymbol\mu=(1,1,1)$ and \protect\subref{subfig:pipeline_3_2_5} $\boldsymbol\mu=(3,2,5)$.
Note that the scale is not the same in both cases.
The first case gives rise to a perfect linear approximation, which is the expected behavior.
The second one displays a lack of smoothness at the interfaces of the patches.
This is due to the fact that we have $C^{1}$ continuity inside each patch while only continuity is guaranteed at the interfaces.}
\label{fig:examples_pipeline}
\end{figure}

Finally, in figure \ref{fig:examples_pipeline}, we show the solution on the whole domain for different values of the parameters.
Note that the solution is not completely smooth at the interfaces of the patches.
This comes from the fact that we have $C^{1}$ continuity in each of the patch while we only have $C^{0}$ continuity at the interfaces.
Methods exist to obtain more regularity at the interfaces (see e.g. \cite{cottrell2009isogeometric}).

\subsection{Geometrical Parameters for Heat Conduction in a Cylinder}\label{subsec:camembert}

\begin{figure}[!ht]
\centering
\subfloat[]{\label{subfig:camembert_original_domain}\includegraphics[clip=true,width=0.2\linewidth]{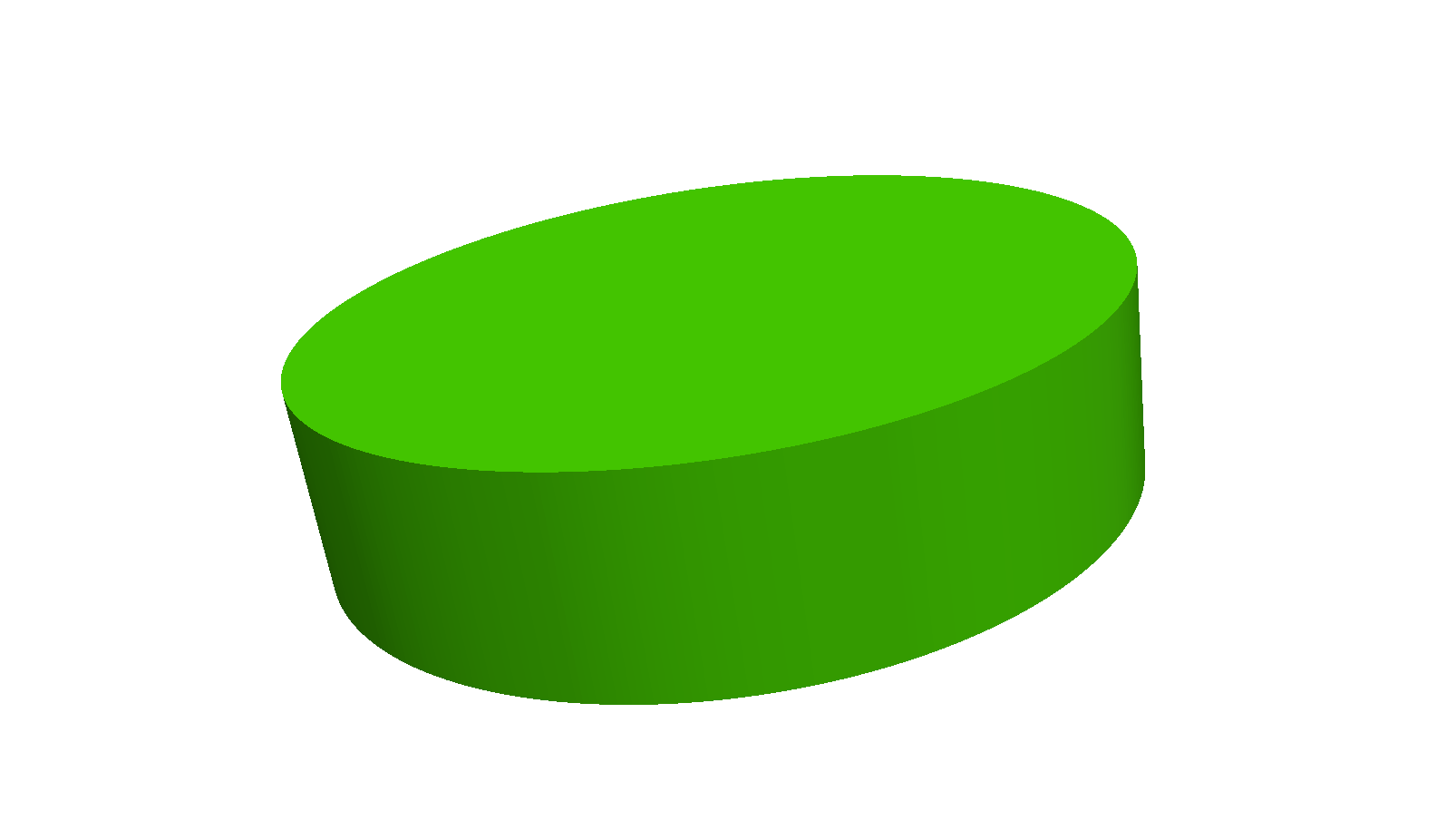}}\
\subfloat[]{\label{subfig:mp_geo_camembert_1_3_4}\includegraphics[clip=true,width=0.2\linewidth]{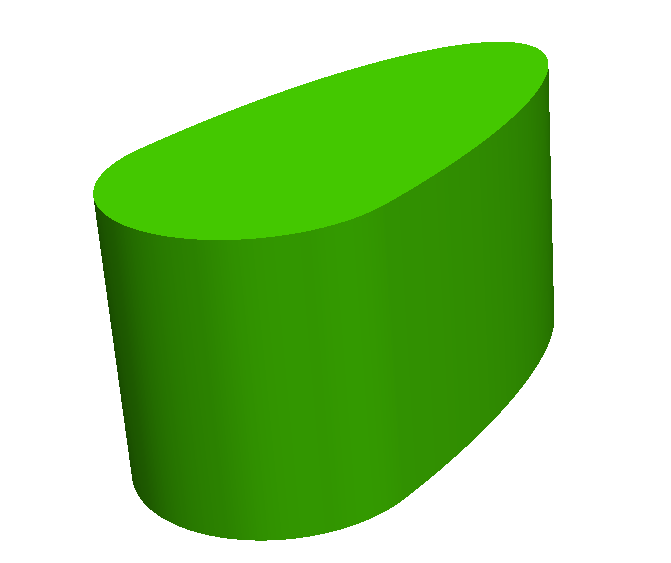}}\
\subfloat[]{\label{subfig:mp_geo_camembert_3_5_2}\includegraphics[clip=true,width=0.2\linewidth]{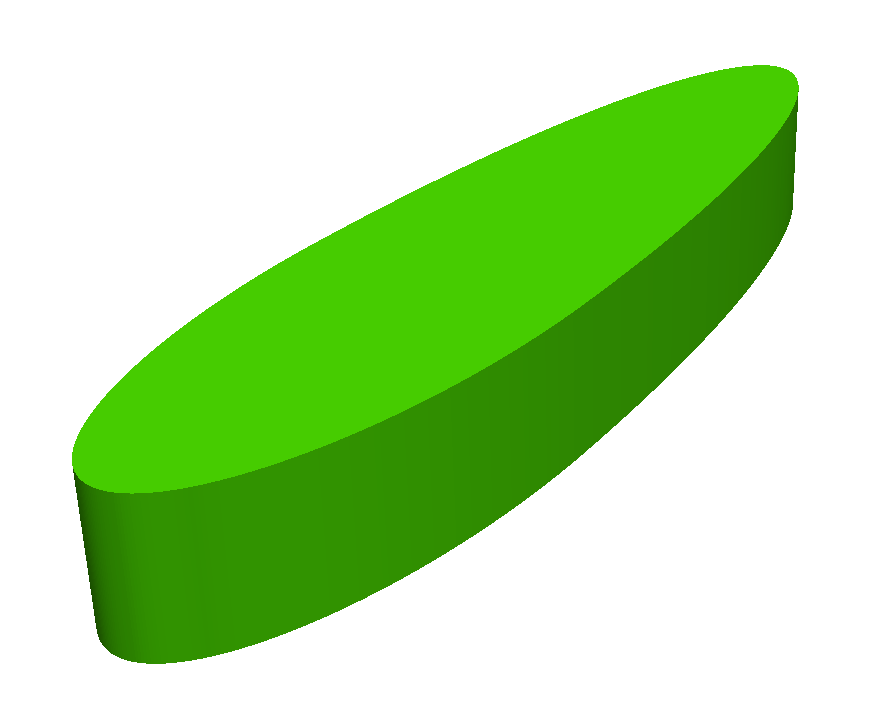}}\
\subfloat[]{\label{subfig:mp_geo_camembert_1_4_1}\includegraphics[clip=true,width=0.2\linewidth]{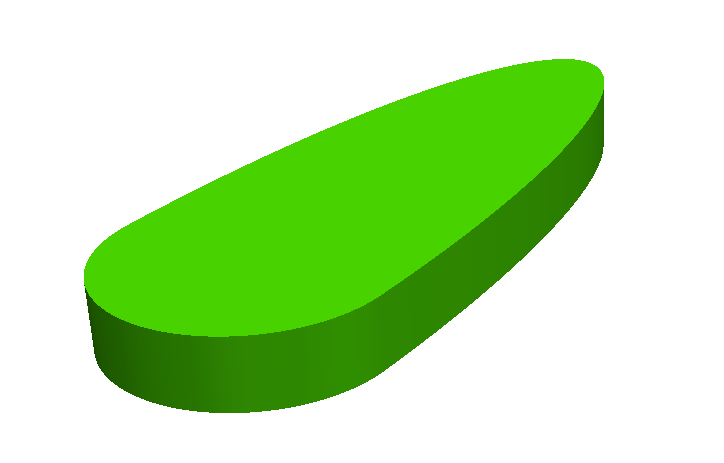}}
\caption{Computational domain for the cylinder test case for different values of the parameters.
The original domain is presented in \protect\subref{subfig:camembert_original_domain}.
Four patches were necessary to build the structure, one for each quarter of the cylinder.
The affine transformation from Section \ref{subsec:camembert} were considered for \protect\subref{subfig:mp_geo_camembert_1_3_4} $\boldsymbol\mu=(1,3,4)$, \protect\subref{subfig:mp_geo_camembert_3_5_2} $\boldsymbol\mu=(3,5,2)$ and \protect\subref{subfig:mp_geo_camembert_1_4_1} $\boldsymbol\mu=(1,4,1)$.}
\label{fig:camembert_domain_transformation}
\end{figure}

We present here the case of a parameter dependent geometry.
The reference domain is a cylinder of radius $2$ and height $1$ oriented in the $z$ direction, i.e. $\Omega(\boldsymbol\mu_{\text{ref}}):=\left\{(r,\theta,z)\ \middle|\ r\in[0,2], \theta\in[0,2\pi],z\in[0,1]\right\}$ where $(r,\theta,z)$ denote the cylindrical coordinates in $\mathbb{R}^{3}$.
The reference cylinder is depicted in Figure \ref{subfig:camembert_original_domain}.
To build it, four patches were necessary.

The transformations under consideration are scaling with respect to the $y$ and $z$ axis.
More precisely, three parameters $\boldsymbol\mu=\left(\mu_{1},\mu_{2},\mu_{3}\right)\in[1,5]^{3}$ are considered, where $\mu_{1}$ scales the portion of the domain for which $y>0$, $\mu_{2}$ the one where $y<0$, and $\mu_{3}$ scales in the $z$ direction.
In other words, the transformation in the part of the domain where $y>0$ is given by
\begin{equation*}
C^{k}=\left(
\begin{array}{cc}
0\\
0\\
0
\end{array}
\right),\qquad
G^{k}=\left(
\begin{array}{ccc}
1 & 0 & 0\\
0 & \mu_{1} & 0\\
0 & 0 & \mu_{3}
\end{array}
\right),
\end{equation*}
while in the region $y<0$ it reads
\begin{equation*}
C^{k}=\left(
\begin{array}{cc}
0\\
0\\
0
\end{array}
\right),\qquad
G^{k}=\left(
\begin{array}{ccc}
1 & 0 & 0\\
0 & \mu_{2} & 0\\
0 & 0 & \mu_{3}
\end{array}
\right).
\end{equation*}
In Figure \ref{fig:camembert_domain_transformation}, we present the domain after application of the affine transformation for different values of the parameters.

%

Considering the conductivity coefficient in \eqref{equ:elliptic_numerics}, we have $a(\boldsymbol\mu,x)=1$.
In order to describe the boundary conditions, let us denote by $\Gamma_{\text{bot}}$, $\Gamma_{\text{top}}$ and $\Gamma_{\text{curve}}$ the bottom, top and curvy boundaries, respectively.
We impose homogeneous Dirichlet boundary conditions on $\Gamma_{D}:=\Gamma_{\text{curve}}$ and unitary Neumann conditions on $\Gamma_{N}:=\Gamma_{\text{top}}\cup\Gamma_{\text{bot}}$.
Finally the right-hand side function is $f=10\chi_{B}$, where $B$ is the ball of radius $0.2$ centered at $(0,0,0.5)$.

\begin{figure}[!ht]
\centering
\includegraphics[clip=true,width=0.35\linewidth]{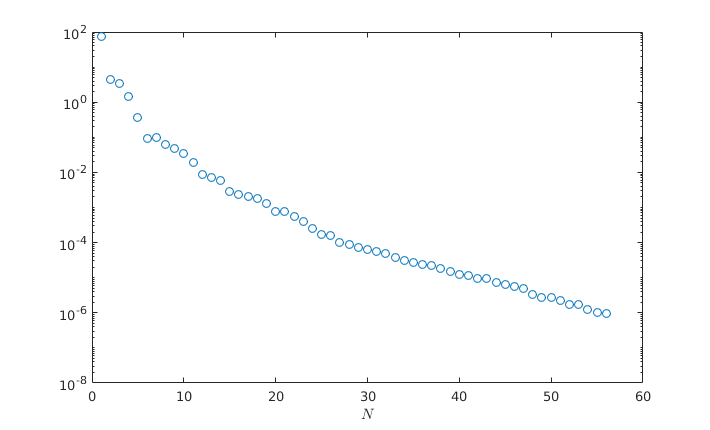}
\caption{Convergence of the greedy algorithm (see Section \ref{subsec:greedy}) for the cylinder test case of Section \ref{subsec:camembert}.}
\label{fig:conv_greedy_camembert}
\end{figure}

Turning to the computational costs, the number of degrees of freedom for the NURBS approximation is $\mathcal{N}=3240$ and the one of the RB is $N=57$.
The whole offline procedure took $30$ minutes while the average RB evaluation takes $5\cdot10^{-4}$ seconds.
In Figure \ref{fig:conv_greedy_camembert}, we show the convergence of the greedy algorithm for this case.

\begin{figure}[!ht]
\centering
\subfloat[]{\label{subfig:camembert_1_4_1}\includegraphics[clip=true,width=0.32\linewidth]{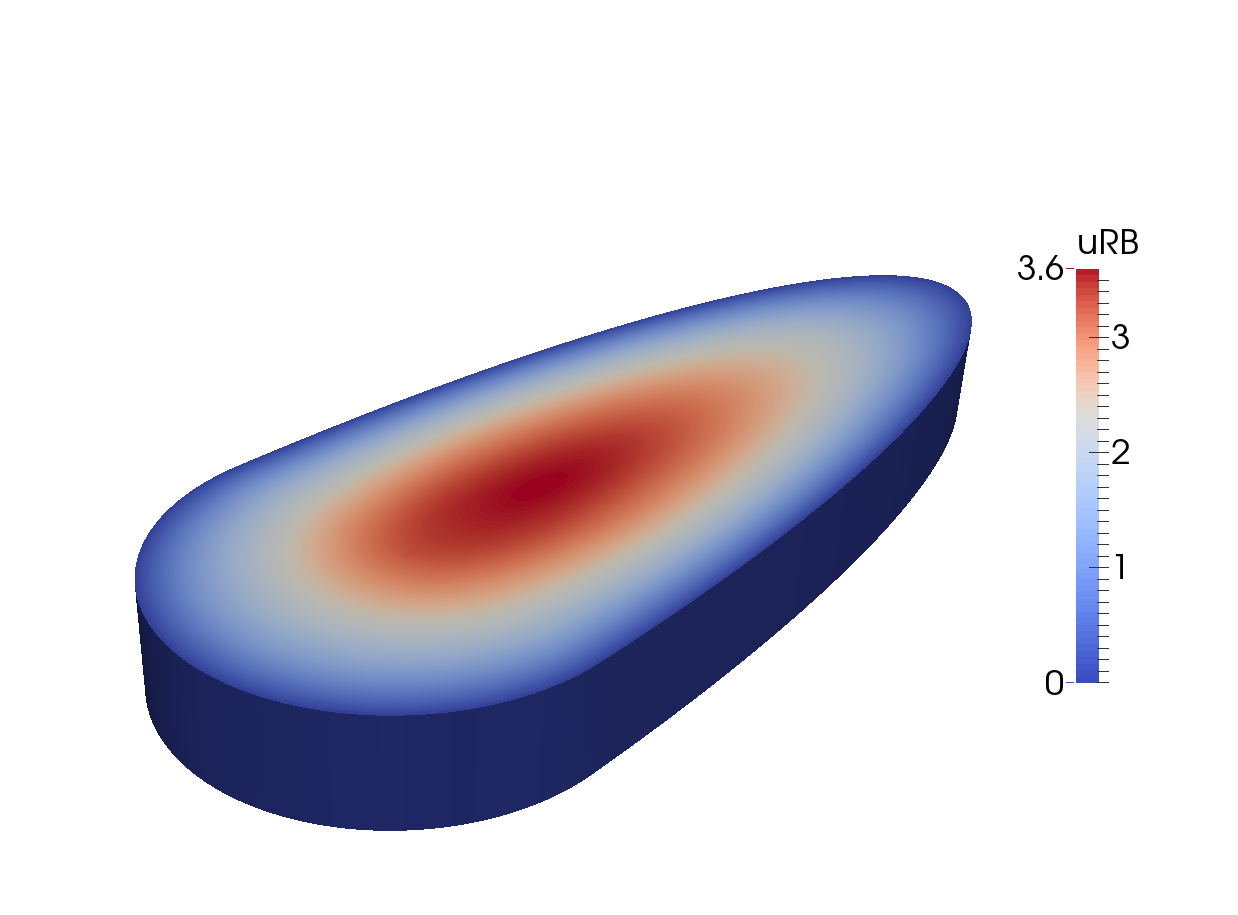}}\qquad
\subfloat[]{\label{subfig:camembert_1_4_1_cut}\includegraphics[clip=true,width=0.32\linewidth]{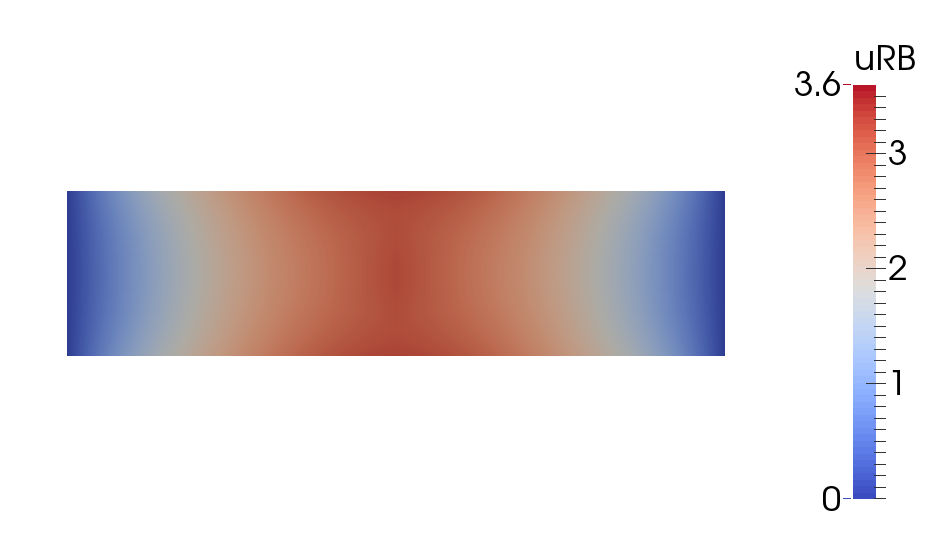}}\\
\subfloat[]{\label{subfig:camembert_1_1_4}\includegraphics[clip=true,width=0.32\linewidth]{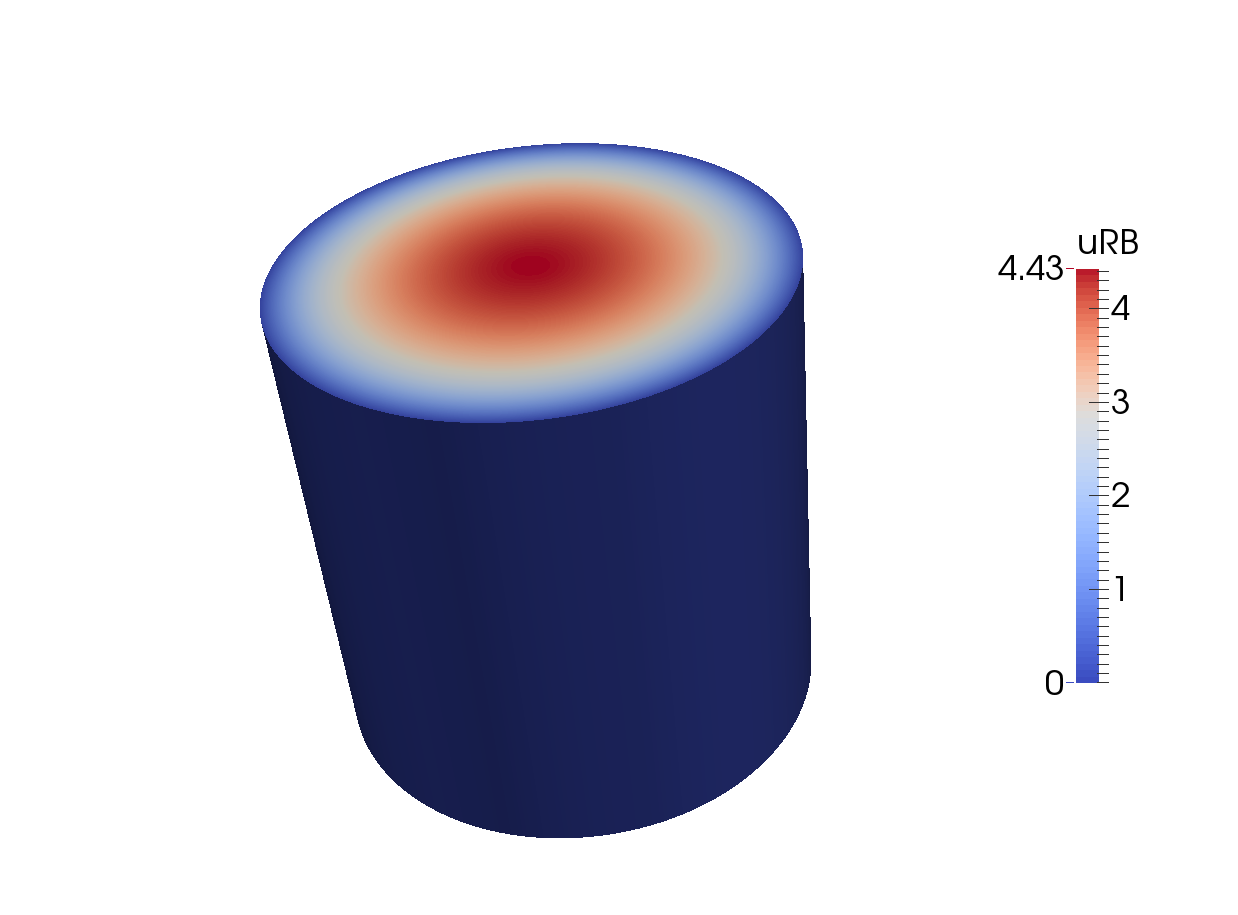}}\qquad
\subfloat[]{\label{subfig:camembert_1_1_4_cut}\includegraphics[clip=true,width=0.32\linewidth]{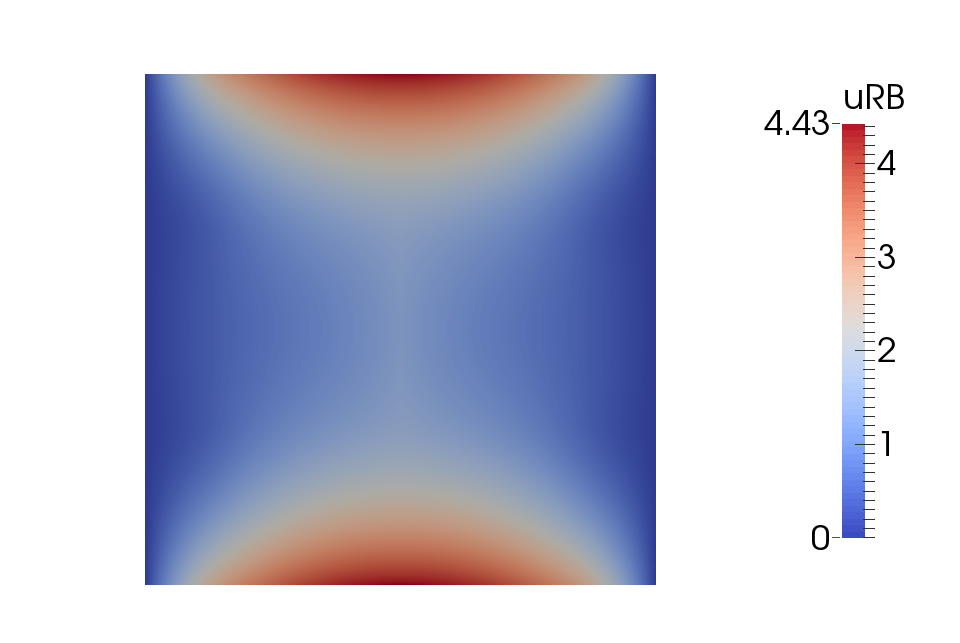}}
\caption{Reduced basis approximation of the cylinder test case for \protect\subref{subfig:camembert_1_4_1} $\boldsymbol\mu=(1,4,1)$ and \protect\subref{subfig:camembert_1_1_4} $\boldsymbol\mu=(1,1,4)$.
The pictures \protect\subref{subfig:camembert_1_4_1_cut} and \protect\subref{subfig:camembert_1_1_4_cut} represent the value of the field on the plane $\left\{(x,y,z)\in\Omega\ \middle|\ y=0\right\}$ for the values considered in \protect\subref{subfig:camembert_1_4_1} and \protect\subref{subfig:camembert_1_1_4}, respectively.
Note that different scales have been used for the different values of the parameters.}
\label{fig:examples_camembert}
\end{figure}

The solution for several values is depicted in Figure \ref{fig:examples_camembert}.
We present the RB approximation on the whole domain as well as its evaluation on the plane $\left\{(x,y,z)\in\Omega\ \middle|\ y=0\right\}$.

\newpage
\bibliographystyle{plain}
\bibliography{Part5_biblio}

\end{document}